\documentclass{aims} 
\usepackage{amsmath}
  \usepackage{paralist}
  \usepackage{graphics} 
  \usepackage{epsfig} 
\usepackage{graphicx}  \usepackage{epstopdf} 
 \usepackage[colorlinks=true]{hyperref}
\hypersetup{urlcolor=blue, citecolor=red}
\allowdisplaybreaks

\def\currentvolume{X}
 \def\currentissue{X}
  \def\currentyear{200X}
   \def\currentmonth{XX}
    \def\ppages{X--XX}

\usepackage{subfigure}
\usepackage{tikz-cd}

\newcommand{\R}{\mathbb{R}}
\newcommand{\Z}{\mathbb{Z}}
\newcommand{\C}{\mathbb{C}}
\newcommand{\N}{\mathbb{N}}

\newcommand{\rmi}{\mathrm{i}}
\newcommand{\rme}{\mathrm{e}}
\newcommand{\rmO}{\mathrm{O}}
\newcommand{\rmId}{\mathrm{Id}}

\newtheorem{innercustomprop}{Proposition}
\newenvironment{customprop}[1]
  {\renewcommand\theinnercustomprop{#1}\innercustomprop}
  {\endinnercustomprop}

  \newtheorem{innercustomthm}{Theorem}
\newenvironment{customthm}[1]
  {\renewcommand\theinnercustomthm{#1}\innercustomthm}
  {\endinnercustomthm}

\newtheorem{theorem}{Theorem}[section]

\newtheorem{lemma}[theorem]{Lemma}
\newtheorem{proposition}{Proposition}

\theoremstyle{definition}

\newtheorem{remark}{Remark}

\newcommand{\ep}{\varepsilon}

\newtheorem{hypothesis}[theorem]{Hypothesis}
\newtheorem*{hypothesis*}{Hypothesis}

\title[Rotating spirals in oscillatory media with nonlocal interactions]
{Rotating spirals in oscillatory media with nonlocal interactions and their normal form} 

\author[Gabriela Jaramillo]{}

 
 \email{gabriela@math.uh.edu}

\thanks{This work is supported by NSF grant DMS-1911742}

\thanks{$^*$Corresponding author: Gabriela Jaramillo}

\begin{document}
\maketitle

\centerline{\scshape Gabriela Jaramillo$^*$}
\medskip
{\footnotesize
 \centerline{Department of Mathematics, University of Houston}
   \centerline{3551 Cullen Blvd., Houston, TX 77204-3008}
   \centerline{USA}
} 



\bigskip

 \centerline{(Communicated by the associate editor name)}


\begin{abstract}
Biological and physical systems that can be classified as oscillatory media 
give rise to interesting phenomena
like target patterns  and spiral waves.
The existence of these structures has been proven in the
case of systems with local diffusive interactions.
In this paper the more general case of oscillatory media with nonlocal coupling is considered.
We model these systems using evolution equations where the nonlocal interactions are
expressed via a diffusive convolution kernel, and prove the existence of
rotating wave solutions for these systems.
Since the nonlocal nature of the equations precludes the use of standard
techniques from spatial dynamics, the method we use relies instead on a combination of
a multiple-scales analysis and a construction similar to Lyapunov-Schmidt. 
This approach then allows us to derive a normal form, or reduced equation, that
captures the leading order behavior of these solutions.
\end{abstract}


\section{Introduction}

In this paper we consider nonlocal evolution equations of the form 
\begin{equation}\label{e:main1}
 \partial_t U =  K \ast U +F(U; \lambda) \quad x \in \R^2,\quad U \in \R^2, \quad \lambda \in \R,
 \end{equation}
where
$K \ast$ represents a component-wise 
convolution operator, and the reaction terms $F(U;\lambda)$ 
undergo a Hopf bifurcation as the parameter 
$\lambda $ is varied. 
This type of nonlocal equations can describe for example oscillating chemical reactions
with components that evolve at different time scales, and where 
 the nonlocal interaction is the result of adiabatically eliminating the fast variable  \cite{shima2004, sheintuch1997, nicola2002}.
Other examples that can be described by similar equations,
but perhaps with different nonlinearities 
include electrochemical systems, \cite{christoph2002, garcia2008}, 
neural field models, \cite{coombes2005, bressloff2011, pinto2001},
and population models where individuals disperse, or move,
 in a nonlocal manner. See for example \cite{kolmogorov1937, coville2007,hutson2003, bates2007} 
 or \cite{andreu2010} and the references therein. 

More generally, equation \eqref{e:main1} is an abstract model for 
oscillatory media. As such it describes any system that is  composed of small oscillating elements 
that are coupled together in a nonlocal  fashion. This coupling is described by the convolution kernel, $K$,
 and the assumptions on this kernel depend on the application. 
Here we focus on general diffusive kernels, $K(x) = -1 + K_1(|x|)$, where  $K_1(|x|)$ is
 exponentially decaying, can take on negative values, and has a finite second moment.
Therefore the convolution operators $K\ast$ describe a process that lies somewhere between anomalous diffusion,
modeled using the fractional Laplacian, 
\[ (- \Delta)^s u(x) = c \int_{\R^2} \frac{f(x) -f(y)}{|x-y|^{2 +2 s} }\;dy, \]
with appropriate constant $c>0$, and regular diffusion, modeled using the standard Laplacian,
 \[ \Delta u(x,y) \sim \frac{1}{h^2}\Big (u(x + h,y) + u(x,y+h) - 4 u(x,y) + u(x-h,y) + u(x,y-h) \Big),\]
 with $(x,y) \in \R^2$, which samples only nearest neighbors.

\begin{figure}
\subfigure[]{\includegraphics[width=.37\textwidth]{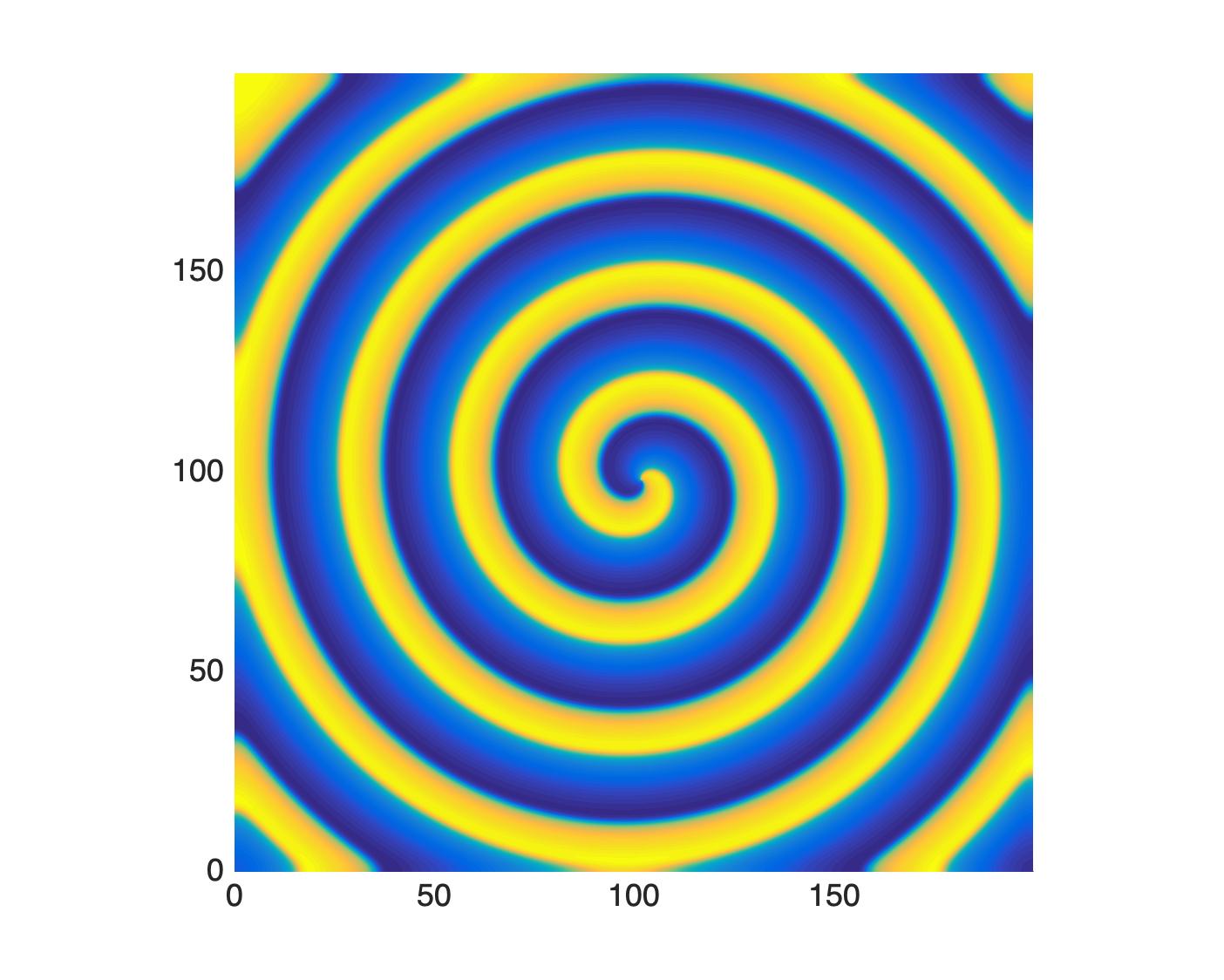}}
\subfigure[]{\includegraphics[width=.4\textwidth]{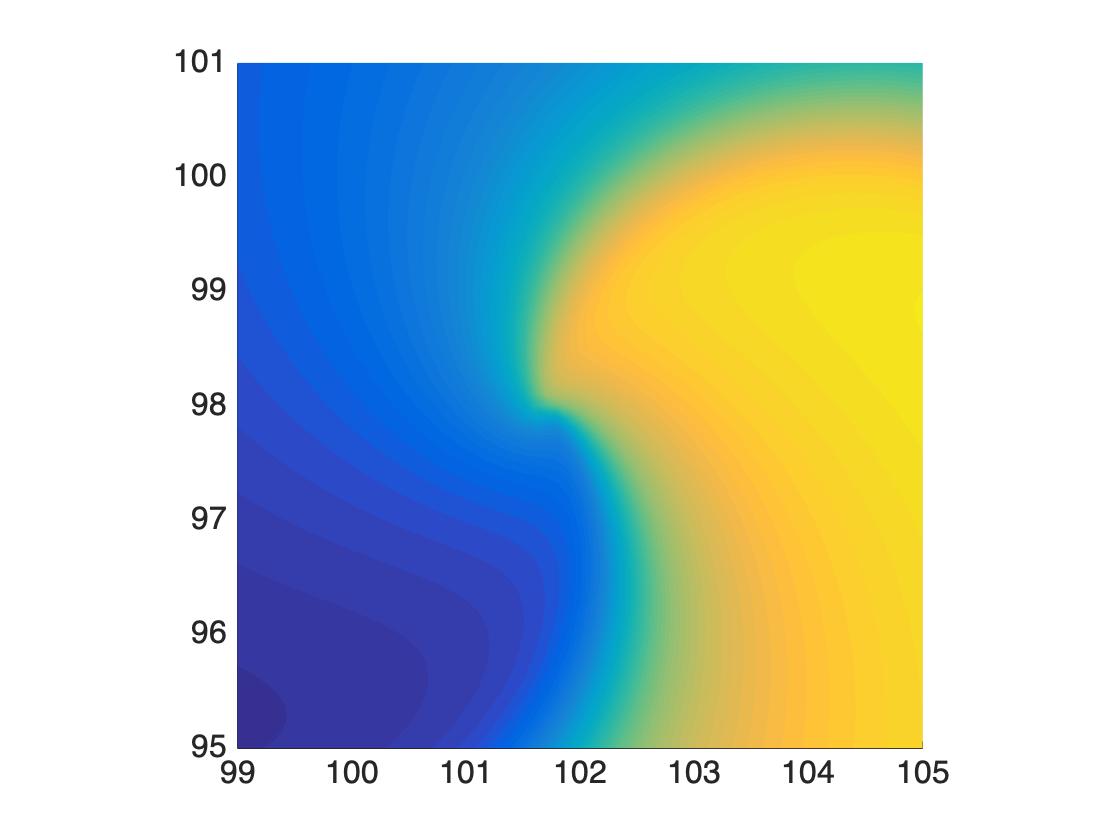}}
\subfigure[]{\includegraphics[width=.37\textwidth]{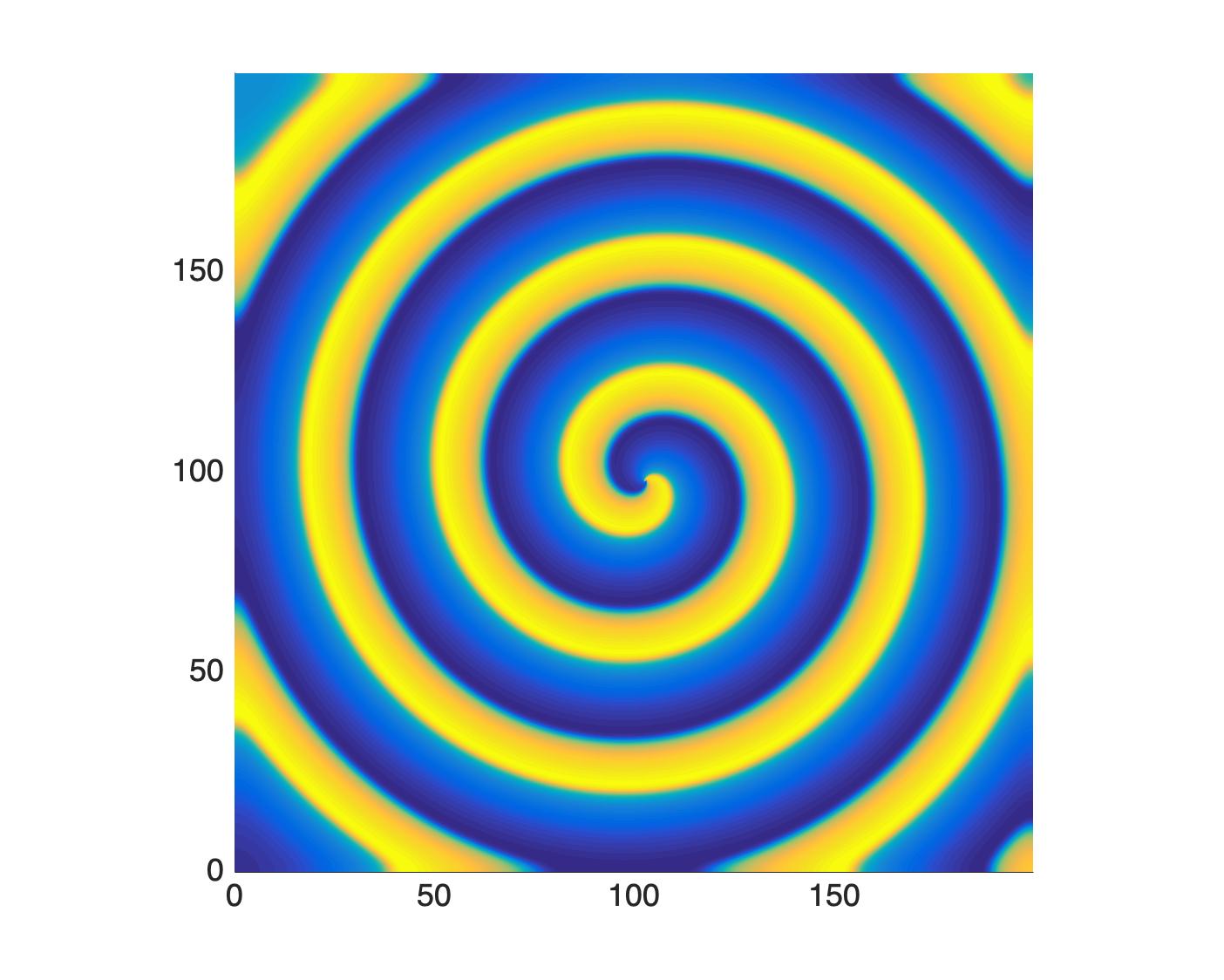}}
\subfigure[]{\includegraphics[width=.4\textwidth]{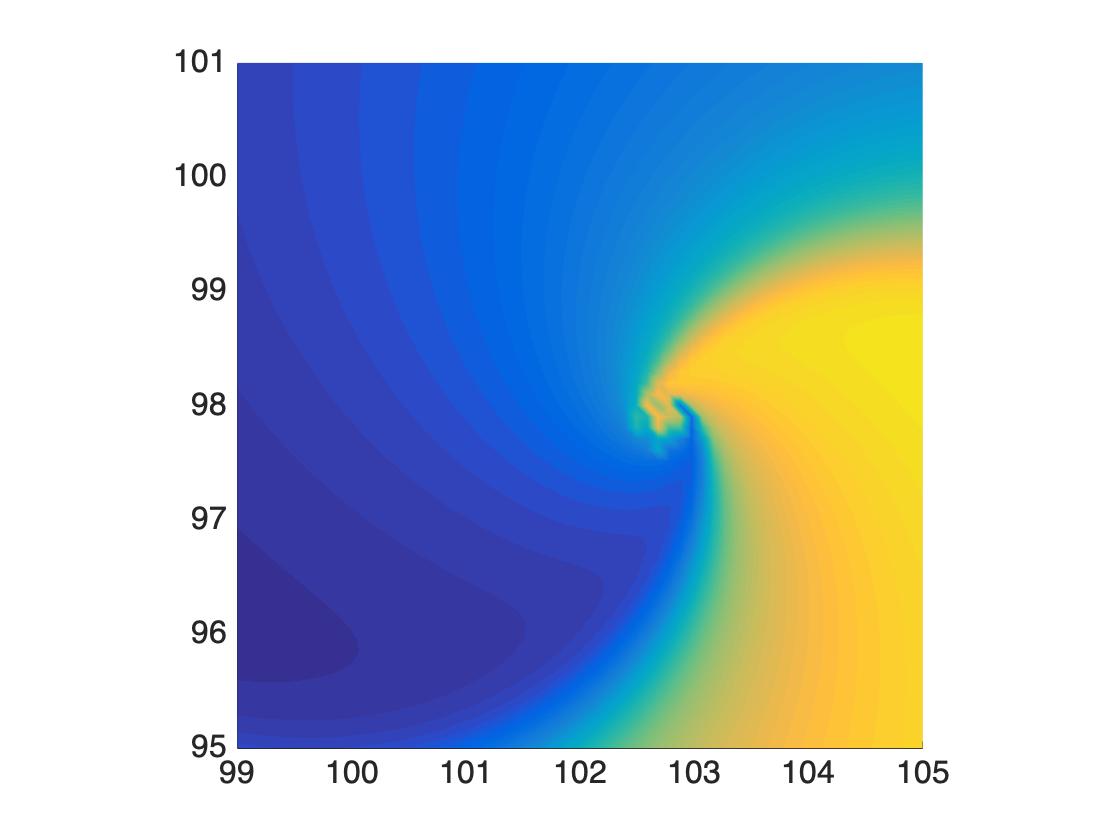}}
\caption{ Spirals obtained by numerical simulations of a FitzHugh-Nagumo system (see Section 6) 
 with $\delta =0.2, \beta =1, \tau =0.1$, and
$K \ast u = \frac{\sigma}{D} [ - u + K_0(|x|/ \sqrt{D}) \ast u]$.
Here $K_0$ is the Modified Bessel function of the second kind.
 In figures a) and b) the choice of $\sigma =5 $ and $D = 0.5$,
 results in a spiral pattern. In contrast, in figures c) and d) the choice of $\sigma=5$ and $D = 1$
 results in a spiral chimera pattern (incoherent core). 
The figures on the right zoom in into the core of the spirals appearing on the left.}
\label{Fig:Spirals}
\end{figure}

These assumptions imply that the processes considered here 
exhibit a characteristic length scale, 
which may suggest that we 
 instead consider reaction diffusion equations
with a rescaled diffusion constant.
However, this simplification misses the true character of this fast diffusion process and precludes 
one from finding interesting patterns like spiral chimeras \cite{shima2004},
chemical turbulence \cite{garcia2008}, localized structures, \cite{colet2014, colet2014-2},
and other phenomena  \cite{plenge2005, nicola2002, siebert2014, scholl2014, nicola2006}. 
Figure ~\ref{Fig:Spirals} also illustrates how slight variations in the choice of convolution kernel, $K$, 
can lead to different patterns in a FitzHugh-Nagumo system with nonlocal diffusion.
The convolution kernel used to create Figures c) and d) involves
longer nonlocal interactions than the one used in Figures a) and b).  
These longer connections give rise to a novel pattern not seen in standard reaction diffusion systems,
 a spiral chimera, which is a spiral pattern with an incoherent core.

Our interest here is on rotating wave solutions,
which we consider as an entry point for understanding more complicated phenomena
in nonlocal systems, like spiral chimeras.
While it is known that rotating waves, and in particular spiral waves,
appear in a wide variety of examples,
 from chemical reactions \cite{nettesheim1993}, to experiments involving slime mold \cite{siegert1995},
 cardiac tissue \cite{davidenko1992, pertsov1993} and brain tissue \cite{Huang2004},
 their existence has only been shown 
 in systems that involve local diffusive processes,
see  \cite{kopell1981,greenberg1980,greenberg1981,hagan1982,scheel1998}.
 In this article we extend these results to oscillatory media with nonlocal coupling. 
 That is, we prove the existence of rotating wave solutions for
problems described by equation  \eqref{e:main1},  and
 derive a normal form which captures their leading order behavior.
  
 In contrast to previous works which use 
reaction diffusion equations \cite{kopell1981, scheel1998},
  the nonlocal character of our problem prevents us from using  techniques
   from  spatial dynamics and dynamical systems theory.
  In particular, we are not able to use center manifold theory as in \cite{scheel1998} to construct
  a normal form and prove existence of solutions in this way.
 Instead, our work is based on the approach taken in the physics
 literature, where a multiple-scale analysis is used to derive 
 amplitude equations \cite{kuramoto2003, tanaka2003, garcia2008}.
 Our main contribution is to place this method in 
 an appropriate functional analytic setting that allows us to use 
 a reduction similar to  Lyapunov-Schmidt  to rigorously derive the normal form
  and obtain its validity at the same time.
In the process we also show the existence
of general rotating wave solutions.
We point out here that to obtain the existence of specific rotating wave solutions, 
like for example spiral waves, one then has to work with the normal form and
show that such patterns exist as solutions to this reduced equation.
  
 Before describing our methods in more detail, we first state our main assumptions
 and describe our set up. Throughout the paper we use polar coordinates $(r, \vartheta)$ and
we look for  rigidly rotating solutions that satisfy $U(r,\vartheta,t) = U(r, \vartheta +ct)$,
 and that are periodic in their second argument, i.e. $U(r, \theta + 2\pi) = U(r, \theta)$. 
We do not consider meandering solutions.  With this notation the problem reduces to 
  solving the steady state equation
  \begin{equation}\label{e:main_steady}
 0 = \underbrace{ K \ast U - c \partial_\theta U + D_UF(0;0)U}_{L} + \underbrace{ [ F(U;\lambda) - D_UF(0,0)U]}_{\mathcal{N(U;\lambda)}},
 \end{equation}
where the rotational speed, $c \in \R$, is treated as an extra parameter.
We further assume that 
\begin{enumerate}
\item $F$ depends on $U$, but not its derivatives.
\item $F(0;\lambda) =0$ for all $\lambda\in \R$.
\item  $D_UF(0;0) =A  \in \mathbb{C}^{2\times 2}$ has a pair of purely imaginary eigenvalues $\pm \rmi \omega$.
\end{enumerate}
We also make the following assumptions on the convolution kernel $K$, which will allow us to show that the operator $L$ defines a Fredholm operator in appropriate weighted spaces.
\begin{hypothesis*}[H1]
The convolution kernel $K$ has a radially symmetric Fourier symbol 
$\hat{K}(\xi) = \hat{K}(|\xi|)$. As a function of $\rho = |\xi|$, 
the symbol  $\hat{K}(\rho)$ can be extended to a uniformly bounded and analytic function on a strip 
$\Omega = \R \times ( -\xi_0, \xi_0) \subset \mathbb{C} $, for some constant $\xi_0>0$.
\end{hypothesis*}

\begin{hypothesis*}[H2]
The symbol $\hat{K}(\rho)$ is symmetric and has a simple zero, 
which we assume is located at the origin $\rho=0$. 
This zero is of order $\ell =2$ and thus $\hat{K}(\rho)$  has the following 
Taylor expansion near the origin:
\[ \hat{K}(\rho) \sim - \alpha \rho^{2} + \rmO( \rho^{4}), \quad \alpha>0.\]
\end{hypothesis*}

As shown in Section \ref{s:convolution}, Lemma \ref{l:decomposition}, with the above hypotheses
the symbol $\hat{K}$ can be written as
\[ \hat{K}(\xi) = \hat{M}_L(|\xi|) \frac{|\xi|^2}{1+ |\xi|^2} =\frac{|\xi|^2}{1+ |\xi|^2} \hat{M}_R(|\xi|),\]
where the functions $\hat{M}_{L/R}(\rho)$, are analytic and uniformly bounded. 
With this information, we make one further assumption.

\begin{hypothesis*}[H3]
The symbols, $\hat{M}_{L/R}$, 
satisfy  $$\hat{M}_L(\rho) = \hat{M}_R(\rho) = \hat{M}(\rho) = c( 1+ \hat{G}(\rho)),$$ with 
$c \in \R$, and $\hat{G}(\rho) \sim \rmO(1/\rho)$ as $\rho \to \infty$.
\end{hypothesis*}

The following are examples of convolution operators that satisfying the above hypotheses:
\begin{enumerate}[1)]
\item $ K \ast u = - \frac{1}{2}  u + \exp(-|x|^2) \ast u$, which has symbol 
$\hat{K}(\xi) = -\frac{1}{2} + \frac{1}{2} \exp(-|\xi|^2/4)$.
\item $K \ast u = \frac{\sigma}{D} [ - u + K_0(|x|/ \sqrt{D}) \ast u]$,
 where $\sigma, D >0$ and $K_0$ represents the Modified Bessel function of the second kind. 
 In this case $K$ has symbol
$\hat{K}(\xi) = \sigma \left[ \frac{-|\xi|^2}{ 1 + D |\xi|^2} \right]$.
\end{enumerate}


\subsection{The Method}
We start this subsection by giving a short overview of the multiple-scales method
(following  \cite{kuramoto2003}),
and then move on to describe our approach.

 For the case of oscillatory reaction diffusion equations in finite domains,
the multiple-scales method is based on the following analysis.
We know that as the parameter $\lambda$ crosses its critical value, the homogenous steady state
 undergoes a Hopf bifurcation. Since we are considering the case of 
 a finite domain, linearizing the equation about this uniform state results 
 in an operator with point spectrum, $\nu_k$, and corresponding eigenmodes, $\rme^{\rmi kx}$ 
 (assuming periodic boundary conditions).
  In particular, we have a pair of complex eigenvalues, $\nu_0 = \lambda \pm \rmi \omega$, associated with the zero mode,
 that cross the imaginary axis when $\lambda >0$.  
 If the system size is small, the eigenvalues are well separated
 and one can study the dynamics of the emerging oscillations by
 keeping track of only the critical mode. 
 This results in the Stuart-Landau equation (an o.d.e.) as a normal form. 
 
 However, as the system size increases, the eigenvalues are no longer well separated. 
 Nonuniform perturbations with small wavenumbers  become important for determining the behavior of solutions that emerge from the bifurcation. 
In some cases it might be possible to describe these solutions
by considering the dynamics of only a finite number of modes, i.e modes with  $|k|<<\delta$, 
 for some small $\delta$, 
and as a result  the normal form becomes a system of o.d.e.

In practice, it is often the case that all modes need to be taken into account,
and instead of separating eigenmodes into critical and rapidly decaying ones,
one separates scales. That is, one establishes 
  fast variables, $x, t,$ and slow  variables, $ X = \ep x, T= \ep^2t$, and approximates the solution 
  as a regular expansion of the form
\[ \tilde{U}(t, X,T; \ep) = \ep U_1(t,X,T) + \ep^2 U_2(t,X,T)+ \ep^3 U_3(t,X,T) +\cdots\]
with 
\[ U_1(t,X,T) = W(X,T) \rme^{\rmi \omega t} + \overline{W}(X,T) \rme^{-\rmi \omega t}.\]
Inserting the approximation $\tilde{U}$ into the reaction diffusion system and collecting terms of equal order in $\ep$ 
results in a hierarchy of equations. 
 Then, applying the Fredholm alternative to all  these equations leads
to solvability conditions, which at order $\rmO(\ep^3)$
  result in an equation for the amplitude, $W(X,T)$.
  This equation is the complex Ginzburg-Landau equation,
   which is now a p.d.e.  and represents the normal form, or amplitude equation, associated with problem. 
  
 The above approach is formal, mainly because this process does not guarantee that 
 the  expansion $\tilde{U} = \ep U_1 + \ep^2 U_2 +\cdots$ converges. 
 In other words, more work needs to be done to show that the complex Ginzburg-Landau equation
  provides valid approximations, $\tilde{U} $,
 for the reaction diffusion system. 
 That is, given a solution, $W$, to the complex Ginzburg-Landau equation 
 that is valid in some time interval $[0,\tilde{T}]$,
 one needs to show that there is a unique solution, $U$, to the original reaction diffusion system
 such that $\|U - U_1\|_{ \mathcal{X}} <\ep^2$ for all $t \in [0,\tilde{T}/\ep^2]$ in some appropriate space $\mathcal{X}$, and where $U_1$ is given as above. 
 Work in this direction is extensive, some of which is
 covered in the references 
 \cite{schneider1992, schneider1995, schneider1996, schneider1998, vanharten1991, kuehn2018}.
Roughly the main idea in these proofs is to show that the residue $R = U - U_1$ remains of order $\rmO(\ep^2)$ on 
the same time interval, $[0,\tilde{T}/\ep^2]$.
For an alternative approach for deriving similar amplitude equations see also \cite{roberts2015}.

Our approach here is to use a similar multiple-scale analysis to derive our normal form,
but because the equation we work with is the steady state equation \eqref{e:main_steady},
 we don't need to look at the evolution of the residue, $R$, to prove the validity of our solutions. 
 Instead, we use a construction similar to a Lyapunov-Schmidt reduction to:
 \begin{enumerate}
 \item[i)] prove the existence of rotating wave solutions to equation \eqref{e:main_steady}, and
 \item[ii)] at the same time obtain the normal form which provides the first order approximation to these solution.
 \end{enumerate}
  We will explain this approach in the remainder of this section.

{\bf Multiple-Scales:} We first separate scales and distinguish between fast variables $r, t$ and slow variables $R= \ep r, T = \ep^2 t$.
Notice that because we are interested in rotating wave solutions, the time scales can be
written directly into the ansatz. That is, we can write
$$U(r,\theta) =U(r,\vartheta +c t) =  U(r, \vartheta + c^*t + \ep^2\mu t ),$$
 where we set  $c = c^* + \ep^2 \mu$. The value of $\mu$ is left as a free parameter and 
the value of $c^*$ is chosen so that given a nonzero integer $n_0$ we have that $ \rmi c^*n_0 = \rmi \omega$, 
i.e. the eigenvalue of the matrix $A = DF(0;0)$. 
Since perturbations with small wavenumbers are the most unstable, 
we assume that the dynamics of the solution is dictated by the slow scales,
 and thus we can write our small amplitude solution as
\begin{equation}\label{e:ansatzU}
 U (r,\theta,R; \ep,\mu) = \ep U_1(\theta,R;\ep,\mu) + \ep^2 U_2(\theta, R;\ep,\mu)  + \ep^3 U_3(r,\theta;\ep,\mu),\end{equation}
with 
\begin{equation}\label{e:ansatzU1}
U_1(\theta,R;\ep,\mu)= W_1 w(R;\ep,\mu) \rme^{\rmi n_0 \theta} + \overline{W}_1 \overline{w}(R;\ep,\mu) \rme^{-\rmi n_0 \theta }. 
\end{equation}
Here the vectors $W_1,\overline{W}_1$ are the eigenvectors of the matrix $A = DF(0;0)$.
Notice also, that the integer $n_0 \in \Z$ is associated with the number of ``arms"  of the rotating wave.
The fact that we can pick a value of $c^*$ such that the relation $ \rmi c^*n_0 = \rmi \omega$ holds,
 implies that rotating wave solutions with any number of ``arms" exists, although they may not be stable.

Inserting the ansatz \eqref{e:ansatzU} into equation \eqref{e:main_steady} 
and separating the different powers of $\ep$ leads to a hierarchy of equations.

 At $\rmO(\ep)$:
 \[ c^*\partial_\theta U_1 - A U_1 =0.\]
 
 At $\rmO(\ep^2)$:
 \[ c^*\partial_\theta U_2 - A U_2 =  MU_1U_1,\] 
 where $M$ represents the Hessian matrix of the reaction term (see Section \ref{s:normalform} for more details).

 And at higher orders:
 \begin{align*}
  c^*\partial_\theta U_3 - K \ast U_3 - A U_3  = & - \ep^2 \mu  \partial_\theta (U_1+  \ep U_2 + \ep^2  U_3) +   \tilde{K}_\ep \ast ( U_1 + \ep U_2 )  \\
  & + \frac{1}{\ep^3} \left[ \mathcal{N}(\tilde{U} ;\lambda) - \ep^2 MU_1U_1\right],
\end{align*}
where due to the separations of scales there is a distinction between $\tilde{K}_\ep$ and $K$. The relation between these two kernels
 is given in terms of their Fourier symbols,
$ \hat{\tilde{K}}_\ep ( P) = \ep^2 \hat{K}( \ep P)$, and more details can be found in Subsection \ref{ss:ms}.
Notice also that our choice of eigenvector $W_1$ and parameters $c^*n_0 = \omega$ 
means that  $U_1$ solves the first equation by design,
and that one can use the second equation to find the function $U_2$ in terms of $U_1$.

{\bf Lyapunov-Schmidt:} The key to our method lies in solving the last equation, which gathers all higher order terms. 
The main point here is that if the left hand side of this equation describes a well behaved operator, 
then one can apply a construction similar to a Lyapunov-Schmidt reduction. More precisely,
in Section \ref{s:normalform} we prove that there is  
projection, $P: X \longrightarrow X_\parallel$, such that the linear operator, 
$$LU = c^*\partial_\theta U - K \ast U - A U ,$$
 can be decomposed into an invertible part $L_\perp: X_\perp \longrightarrow Y$, and a bounded operator $L_\parallel: X_\parallel \longrightarrow Y$. 
 This allows us to split the last equation into two systems,
 \begin{align}\label{e:reduced_intro}
 L_\parallel U_3  = &  P \left[- \ep^2\mu  \partial_\theta (U_1+  \ep U_2 + \ep^2  U_3) +   \tilde{K}_\ep \ast ( U_1 + \ep U_2 ) \right. \\
  \nonumber
  &\qquad \left. + \frac{1}{\ep^3} \left[ \mathcal{N}(\tilde{U} ;\lambda) - \ep^2 MU_1U_1\right] \right],\\[2ex] 
  \label{e:invert_intro}
 L_\perp U_3 = & (\mathrm{I} - P) \left[ - \ep^2 \mu  \partial_\theta (U_1+  \ep U_2 + \ep^2  U_3) +   \tilde{K}_\ep \ast ( U_1 + \ep U_2 ) \right.\\ \nonumber
 & \hspace{10ex} \left.+ \frac{1}{\ep^3} \left[ \mathcal{N}(\tilde{U} ;\lambda) - \ep^2 MU_1U_1\right] \right].
\end{align}
If one further assumes that the term $U_3 \in X_\perp$, 
 then one can use the implicit function theorem to find
 a unique solution $U_3$ to equation \eqref{e:invert_intro}, that depends on $U_1$, $\ep$ and $\mu$.
 Inserting this solution back into \eqref{e:reduced_intro} leads to a reduced equation for $U_1$,
 which after projecting onto the correct space becomes
 an equation for $w$ and defines  the normal form.
 Here we break with tradition and include all terms in this
  reduced equation as part of the definition of the normal form. 
 In Section \ref{s:validity} we will prove that with minimal knowledge of these 
 nonlinearities, mainly that they depend only on $U_1$ and not its derivatives,
  solutions to this normal form exist.
  This result, together with the Lyapunov-Schmidt reduction just described, then implies
  that our normal form is valid.
 We make this statement more precise in the next subsection where we also state our main theorem.

There are two technical difficulties that appear in the above approach. 
One is showing that such a splitting for the linear operator, $L$, is possible, 
and the second one is showing that the nonlinearities are well defined 
functions in the chosen spaces. 
The first difficulty arises because we are setting our equations in the whole plane, $\R^2$. 
Therefore the linearization, $L$,  which behaves very much like a second order elliptic operator,
 has essential spectrum that touches the origin. 
In addition, because of the symmetries present in the problem, 
we have a zero eigenvalue embedded in this essential spectrum. 
As a result, the linearization does not have a closed range
when posed as an operator between 
 regular Sobolev spaces. 
In Section \ref{s:normalform} we will show that if instead we use a special class of weighted Sobolev spaces 
that impose algebraic growth or decay,
we can recover  Fredholm properties for $L$ 
(see  \cite{mcowen1979} for a similar situation in the case of the Laplacian, and \cite{jaramillo2019, jaramillo2018} for the case of a nonlocal equation in one and two dimensions, respectively).
In essence, the extra structure
included in these spaces means that we can remove approximate eigenvalues from our domain
and therefore obtain an operator with closed range, and that
has a finite dimensional kernel (nullspace) and finite dimensional cokernel.

The second difficulty comes from including spiral waves among the rotating wave solutions of interest.
Because these patterns are described by functions which are bounded at infinity,
 we have to consider spaces that allow for a small algebraic growth. 
However, this implies that the chosen spaces are not Banach algebras, and as a result
one is not able to show that the nonlinearities are well defined.
To resolve this issue, we construct spaces that allow us to split our solutions
into the direct sum of algebraically decaying functions plus uniformly bounded functions,
 i.e.  $\mathcal{H} = H^2_\gamma(\R^2) \oplus W^{2,\infty} (\R^2)$, and show
 that the linear operator, as well as the nonlinearities, maintain this structure.

\subsection{ Main Result:}
With assumptions (H1), (H2), (H3), and the approach just described,
we prove that small amplitude rotating wave solutions,
 $U(r,\theta)$, to equation \eqref{e:main_steady} 
 can be approximated by a function  of the form
 \begin{equation}\label{e:ansatz_theorem}
  U_1(r, \theta ; \ep,\mu) =  \left(W_1 w(\ep r;\ep,\mu) \rme^{\rmi n_0 \theta} + \overline{W}_1 \overline{w}(\ep r;\ep,\mu) \rme^{-\rmi n_0 \theta } \right),
  \end{equation}
where $\ep$ and  $\mu$ are small parameters, $c^*n_0 = \omega $,
and the vectors $W_1,\overline{W}_1$ are the eigenvectors of the matrix $A$.
In particular, we show that the function $w$ in this approximation is a solution to
the normal form equation,
 \begin{equation}\label{e:normalform_intro}
 \tilde{K}_\ep \ast w + (\lambda + \rmi (\mu^*+\mu) n_0) w + a|w|^2 w + \rmO(\ep |w|^4 w)=0.
 \end{equation}
Here $\mu^* $ is a small nonzero fixed number, $\lambda $ and $ a $ are appropriate constants
that depend on the system's parameters,
 and as mentioned before, the Fourier symbol for the kernel 
 $\hat{\tilde{K}}_\ep = \frac{1}{\ep^2} \hat{K}(\rho)$, satisfies Hypotheses (H1), (H2), and (H3).

The following proposition, which is proved in Section \ref{s:validity},  then
guarantees the existence of radial solutions, $w(r)$, to this equation,
even without 
explicit knowledge of the higher order terms. In particular, these solutions 
 live in a space $\mathfrak{H}_{\gamma,n} \subset \mathcal{H} = H^2_\gamma(\R^2) \oplus W^{2,\infty} (\R^2)\subset C_B(\R^2)$,
  consisting of functions that can be decomposed into the sum of a bounded function
 and an algebraically localized function.  
 For a more detailed description of the space $\mathfrak{H}_{\gamma,n}$ see Section \ref{s:normalform}.
\begin{proposition}
Given real numbers $ \mu^* \neq 0 $,  $\gamma >0$,
and an integer $n_0$, there exists positive constants $\ep_0, \mu_0$, and a $C^1$ map 
\[
\begin{array}{ c c c c}
\Gamma: & (-\ep_0, \ep_0) \times ( \mu^*-\mu_0, \mu^*+ \mu_0) & \longrightarrow & \mathfrak{H}_{\gamma,n_0}\\
& (\ep,\mu) & \longmapsto & w(R; \ep, \mu)
\end{array}
\]
such that $w(R;\ep,\mu)$ is a solution to equation \eqref{e:normalform_intro}.
\end{proposition}

As a consequence of this proposition, the ansatz \eqref{e:ansatzU}, \eqref{e:ansatzU1}
and the Lyapunov-Schmidt
reduction described in the previous subsection,
 we arrive at our main Theorem.
\begin{theorem}\label{t:main}
Let $\gamma >0 $, $n_0 \in \Z$, and suppose $w(R;\ep,\mu) \in \mathfrak{H}_{\gamma,n_0} $ is a solution to equation \eqref{e:normalform_intro}. 
Then, there exist a unique solution $U(r,\theta)$ of the steady state equation \eqref{e:main_steady} and constants $ C,\ep_*>0$, such that for all $\ep \in (-\ep_*,\ep_*) $ 
the estimate
\[ \| U(r,\theta) - U_1(r,\theta) \|_{C_B} < C\ep^2,\]
with $U_1$ as in \eqref{e:ansatz_theorem}, holds.
\end{theorem}

\begin{remark}
Note that if we let $\theta   = \vartheta +c^*t + \ep^2 (\mu^*+\mu) t$  in  \eqref{e:ansatz_theorem}  ,
 then $U_1$ is an approximation  to a solution of equation \eqref{e:main1}.
 \end{remark}
 \begin{remark}
If the coupling radius is small
compared to the length-scale of the solution/pattern,
then the operator  $K \ast$ can be approximated by the Laplacian and 
the reduced equation \eqref{e:normalform_intro} is the same ordinary differential equation 
obtained in \cite{scheel1998} using center manifold theory.
We suspect that in this case one should be able to adapt the theory presented 
in \cite{faye2018} to show the existence of a center manifold for the nonlocal system \eqref{e:main1}
and thus also obtain the same normal form in this way. 

 \end{remark}
\begin{remark}
Our analysis shows  that a family of rotating wave solutions parametrized by $\ep, \mu$ bifurcates
from the uniform oscillatory state. The parameter $\ep$ relates to the amplitude of these solutions,
whereas the parameter $\mu$ relates to their speed.
The value of $\mu^*$, which appears in the proposition above, is arbitrary;
its only restriction being that it is not equal to zero.
This condition is more a consequence of the method we use to prove the result, i.e. showing
that an operator is invertible, rather than 
some restriction that comes from the problem itself. 
In other words, we can find rotating wave solutions with speeds that are arbitrarily close to $c^*$.
In the case of spiral waves, we suspect that, just like in the case of reaction diffusion equations, 
  there is a particular speed that is selected by the system.
We comment more on this last point in Sections \ref{s:validity} and \ref{s:discussion}.
\end{remark}
\begin{remark}
 As already pointed out, the integer $n_0$ is related to the number of ``arms" of the
rotating wave solution, and the results presented here show that solutions 
with any number of ``arms" exists.
In the case of reaction diffusion systems, it
has been shown, at least formally, that multi-armed spiral waves are
unstable, \cite{hagan1982, greenberg1980,kopell1981}.
We suspect that a similar result holds for spiral waves in oscillatory media with 
nonlocal coupling.
\end{remark}

\subsection{Outline}
The paper is organized as follows: In Section \ref{s:weighted} we introduce a special class of weighted Sobolev spaces and summarize properties of the convolution operator $K$. In Section \ref{s:convolution} we use a Fourier series expansion in the angular variable to obtain a diagonal representation of the linear operator, $L$, which will prove useful when deriving Fredholm properties for this map. The Fredholm properties of $L$ are shown in Subsection \ref{ss:linear_operator}. The rest of Section \ref{s:normalform} is dedicated to deriving a general normal form for rotating solutions using Lyapunov-Schmidt reduction. The validity of the normal form is proved 
 in Section \ref{s:validity},
and in Section \ref{s:example} we determine the reduced equation for a specific example. Finally, we conclude the paper with a discussion in  Section \ref{s:discussion}.

\section{Some Useful Weighted Sobolev Spaces}\label{s:weighted}

In this section we define various weighted Sobolev spaces and their main properties.
These spaces will appear in later sections and are important for showing the
Fredholm properties of the convolution operators considered here.
We will also show that when these weighted spaces are also Hilbert spaces,
 they can be decomposed as a direct sum  which is left invariant under the Fourier Transform.
 This decomposition will allow us to pick the correct projection 
 that is necessary for carrying out the Lyapunov-Schmidt reduction of our main equation.

\subsection{Kondratiev Spaces}

Given $d \in \mathbb{N}$, $s \in \N \cup \{0\}$, $\gamma \in \R$, and $ p \in (1, \infty)$,
we denote Kondratiev spaces by the symbol $M^{s,p}_\gamma(\R^d)$
and define them as the space of locally summable, 
$s$ times weakly differentiable functions with norm
\[ \|u\|^p_{M_\gamma^{s,p}} = \sum_{|\alpha|\leq s } \|  \langle x \rangle^{\gamma + |\alpha|} D^\alpha u \|^p_{L^p}.\]
Here $\alpha$ is a multi-index and we use the notation 
$\langle x \rangle^{\gamma + |\alpha|} = (1 + |x|^2)^{(\gamma + |\alpha|)/2}$. 
Throughout the paper we also identify $M^{0,p}_\gamma(\R^d)$ with $L^p_\gamma(\R^d)$.

To the best of our knowledge Kondratiev spaces where first introduced in
\cite{kondrat1967} to study boundary value problems for elliptic operators. 
Variations of these spaces have also been used
to study exterior Stokes problems, see for example \cite{specovius1986, girault1991, alliot2000}.
More recently, we have used Kondratiev spaces to show existence of target patterns in 
oscillatory media, \cite{jaramillo2015, jaramillo2016, jaramillo2018}
and to study the emergence of striped patterns via inhomogeneities in the 1-d Swift-Hohenberg equation 
\cite{jaramillo2019},
and in the 2-d complex Ginzburg-Landau equation \cite{jaramillo2015-2}.

Notice that depending on the value of $\gamma$, 
these spaces impose a degree of algebraic decay or growth on functions.
In addition, the embeddings $M^{s,p}_\alpha(\R^d) \subset M^{s,p}_\beta(\R^d)$ 
hold provided $\alpha > \beta$, 
and if $s>r$ we also obtain that the space 
$M^{s,p}_\gamma(\R^d) \subset M^{r,p}_\gamma(\R^d)$.
As is the case with regular Sobolev spaces, 
one can also identify the dual $(M^{s,p}_\gamma(\R^d))^*$ 
with the space $M^{-s,q}_{-\gamma}(\R^d)$, where $p,q$ are conjugate exponents.
The notation
\[ (f,g) = \sum_{|\alpha|\leq s} \int_{\R^d} D^\alpha f D^\alpha g\;dx\]
can then be used to denote the pairing between an element 
$f \in M^{s,p}_\gamma(\R^d)$ and an element $g \in M^{-s,q}_{-\gamma}(\R^d)$.
In the particular case when $p =2$, Kondratiev spaces are also Hilbert spaces and have the inner product
\[ \langle f, g \rangle := \sum_{|\alpha|\leq s} \int_{\R^d} D^\alpha f D^\alpha g \cdot ( 1+ |x|^2)^{(\gamma + |\alpha|)} \;dx.\]

The following lemma is the main result of this  subsection.
It establishes that the above Hilbert spaces have a direct sum decomposition
that will prove useful when showing the 
Fredholm properties of the linear part of equation \eqref{e:main_steady}. 
This decomposition holds for all dimension $\R^d$, 
but we restrict our exposition to the two dimensional case, which is relevant here.
In the Lemma we use the symbol $M^{s,p}_{r,\gamma}(\R^2)$
to denote the subset of radially symmetric functions in $M^{s,p}_\gamma(\R^2)$,
 and the overbar to denote the complex conjugate.

\begin{lemma}\label{l:polar}
Given $s\in \N \cup \{0\}$ and $\gamma\in \R$, the space 
$M^{s,2}_\gamma(\R^2)$
can be written as a direct sum decomposition 
\[ M^{s,2}_\gamma(\R^2) = \oplus m^s_{\gamma,n},\]
where $n \in \Z$ and 
 \[ m^s_{\gamma,n}= \{ u \in M^{s,2}_\gamma(\R^2) \mid u(r,\theta) = u_n(r)\rme^{\rmi n \theta} \quad \mbox{and} \quad u_n(r) \in M^{s,2}_{r,\gamma}(\R^2;\C), \quad \bar{u}_n = u_{-n}\}.\]

\end{lemma}

In particular, identifying each point $(x,y) \in \R^2$ with the complex number
$ x + \rmi y = r\rme^{\rmi \theta} $,
the lemma shows that given any $f \in M^{s,2}_\gamma(\R^2)$ one can write
\begin{equation}\label{e:fourierseries}
	f(r \rme^{\rmi \theta} ) = \sum_{n \in \Z} f_n(r) \rme^{\rmi n \theta}\qquad \mbox{with} \qquad f_n(r) = \frac{1}{2\pi} 		\int_0^{2 \pi} f(r \rme^{\rmi \theta}) \rme^{-\rmi n \theta} \;d\theta,
\end{equation}
 where $f_n(r) \in M^{s,p}_{r,\gamma}(\R^2)$.
The proof of this result follows a similar analysis as that of Stein and Weiss 
 \cite[Chapter IV, Section 1]{stein1971}.
We give here a sketch of the proof. 
\begin{proof}
The goal is to show that each $f \in M^{s,2}_\gamma(\R^2)$ 
can be approximated by an element in  $\oplus m^s_{\gamma,n}$.
In particular, the representation \eqref{e:fourierseries} is the desired candidate function.
To start, notice that a short calculation shows that $f_n(r)$
is indeed in $M^{s,p}_{r,\gamma}(\R^2)$.
Next, to prove that
 \[ \| f - \sum_{-N}^N f_n(r) \rme^{\rmi n \theta} \|_{M^{s,2}_\gamma} \rightarrow 0 \quad \mbox{as} \quad n\rightarrow \infty,\]
 we use an equivalent definition for the norm on $M^{s,2}_\gamma(\R^2)$.
 Namely, we consider
 \[ \|f\|_{M^{s,2}_\gamma} = \sum_{|\alpha|\leq s} \left \| \frac{\partial^{|\alpha|}f}{\partial r^{\alpha_1} \partial \theta^{\alpha_2}} \frac{1}{r^{\alpha_2}} \right \|_{L^2_{\gamma+|\alpha|}},\]
 and show that 
 \[ \left\| \frac{\partial^{|\alpha|}f}{\partial r^{\alpha_1} \partial \theta^{\alpha_2}} \frac{1}{r^{\alpha_2}} - \sum_{-N}^N \frac{\partial^{\alpha_1}f_n}{\partial r^{\alpha_1}} \frac{(in)^{\alpha_2}}{r^{\alpha_2}} \rme^{\rmi n \theta} \right \|^2_{L^2_{\gamma+|\alpha|}}\longrightarrow 0 \quad \mbox{as} \quad N \rightarrow \infty.\]

Now, because $f \in M^{s,2}_\gamma(\R^2)$, the integral 
\[ \int_{\R^2} \left |\frac{\partial^{|\alpha|}f}{\partial r^{\alpha_1} \partial \theta^{\alpha_2}} \frac{1}{r^{\alpha_2}} \right|^2 \langle r \rangle^{2(\gamma+ |\alpha|)} \;dx< \infty,\]
so that by Fubini's theorem, the function 
$\frac{\partial^{|\alpha|}f}{\partial r^{\alpha_1} \partial \theta^{\alpha_2}} \frac{1}{r^{\alpha_2}}\langle r \rangle^{(\gamma+ |\alpha|)} $
defines a square integrable function in the variable $\theta$,
and has a Fourier series expansion that converges in the $L^2$ norm for a.e. $r$.
Moreover, one can use integration by parts to show that the
coefficients of this Fourier series are given by 
$\frac{\partial^{\alpha_1}f_n}{\partial r^{\alpha_1}} \frac{(in)^{\alpha_2}}{r^{\alpha_2}} \langle r \rangle^{(\gamma+ |\alpha|)}.$
Using Parseval's theorem, one then obtains that
\[ g_N =  \sum_{-N}^N \left| \frac{\partial^{\alpha_1}f_n}{\partial r^{\alpha_1}} \frac{(in)^{\alpha_2}}{r^{\alpha_2}} \right|^2 \langle r \rangle^{2(\gamma+ |\alpha|)} 
 \longrightarrow \quad
\frac{1}{2\pi} \int_0^{2 \pi} \left| \frac{\partial^{|\alpha|}f}{\partial r^{\alpha_1} \partial \theta^{\alpha_2}} \frac{1}{r^{\alpha_2}} \right|^2 \langle r \rangle^{2(\gamma+ |\alpha|)} \;d\theta\]
for a.e. $r$, which allows us to use  the Dominated convergence theorem and conclude that
\[
\begin{split}
  \lim_{N \to \infty} \int_0^\infty \sum_{-N}^N  \left| \frac{\partial^{\alpha_1}f_n}{\partial r^{\alpha_1}} \frac{(in)^{\alpha_2}}{r^{\alpha_2}} \right|^2 & \langle r \rangle^{2 (\gamma+ |\alpha|)} r \; dr =\\
   \frac{1}{2 \pi}  \int_0^\infty \int_0^{2\pi} &  \left| \frac{\partial^{|\alpha|}f}{\partial r^{\alpha_1} \partial \theta^{\alpha_2}} \frac{1}{r^{\alpha_2}} \right|^2 \;d\theta\langle r \rangle^{2 (\gamma+ |\alpha|)} r \; dr .
 \end{split}\]
 
Next, we have that
\[
\begin{split} 
\int_{\R^2} \left | \frac{\partial^{|\alpha|}f}{\partial r^{\alpha_1} \partial \theta^{\alpha_2}} \frac{1}{r^{\alpha_2}}  - \sum_{-N}^N \frac{\partial^{\alpha_1}f_n}{\partial r^{\alpha_1}} \frac{(in)^{\alpha_2}}{r^{\alpha_2}} \rme^{\rmi n \theta}\right |^2 &  \langle x \rangle^{2( \gamma+|\alpha|)} \;dx  = \\
\int_0^\infty  \left[ \int_0^{2\pi}  \left| \frac{\partial^{|\alpha|}f}{\partial r^{\alpha_1}  \partial \theta^{\alpha_2}} \frac{1}{r^{\alpha_2}} \right|^2 \right. & \left. \langle r \rangle^{2 (\gamma+|\alpha|)}  \;d\theta - g_N  \right]r \; dr,
\end{split}
 \]
which is obtained by expanding the integrand in the left hand side,
and integrating with respect to $\theta$.
The result now follows by combining these last two equalities.

\end{proof}

The next lemma makes explicit the relation between the parameter $\gamma$
and the decay rate of elements in the space $M^{2,2}_\gamma(\R^d)$.
In particular, it shows that given any function $f \in M^{2,2}_\gamma(\R^2)$,
then $f$ decays algebraically for values of $\gamma>-1$ 
and grows algebraically for values of $\gamma<-1$.
This result will prove useful in Section \ref{s:normalform}.
\begin{lemma}\label{l:growth}
Suppose $f \in M^{2,2}_\gamma(\R^d)$ then 
$|f(x)| \leq C \|f \|_{M^{2,2}_\gamma} |x|^{-(\gamma+d/2)}$,
with $C$ a generic constant.
\end{lemma}

A proof of this lemma can be found in the Appendix.


\subsection{Weighted Sobolev Spaces:}
Throughout the paper we will also use the notation 
$W^{s,p}_\gamma(\R^d)$ 
to denote the space of locally summable functions with norm
\[ \|u\|^p_{W_\gamma^{s,p}} = \sum_{|\alpha| \leq s } \| \langle x \rangle^{\gamma} D^\alpha u \|^p_{L^p},\]
where $\alpha$ is again a multi-index. 
When $p=2$ to simplify notation we let 
$H^s_\gamma(\R^d) = W^{s,p}_\gamma(\R^d)$. 
Notice that in contrast to Kondratiev spaces, functions in the space 
$W^{s,p}_\gamma(\R^d)$ do not necessarily gain localization with each derivative.  
The following lemma shows that the spaces $H^s_\gamma(\R^d)$
can also be decomposed as a direct sum.
As before, we use the notation $H^s_{r, \gamma}(\R^2)$ to denote 
those functions in $H^s_\gamma(\R^2)$ that are radially symmetric.
 
\begin{lemma}\label{l:polar2}
Given $s\in \N \cup \{0\}$ and $\gamma\in \R$,
the space $H^s_\gamma(\R^2)$ can be written as a direct sum decomposition 
\[ H^s_\gamma(\R^2) = \oplus h^s_{\gamma,n},\]
where $n \in \Z$ and 
\[ h^s_{\gamma,n}= \{ u \in H^s_\gamma(\R^2) \mid u(r,\theta) = u_n(r)\rme^{\rmi n \theta} \quad \mbox{and} \quad u_n(r) \in H^s_{r,\gamma}(\R^2;\C), \quad \bar{u}_n = u_{-n}\}.\]
\end{lemma}
Since the proof of this result follows a similar argument as 
that of Lemma \ref{l:polar}, we omit the details.

\begin{remark}\label{r:normH}
Notice that the space $ h^s_{\gamma,n}$ comes equipped with the norm
 \[ \| u(r,\theta)\|_{h^s_{\gamma,n}} =  \sum_{| \alpha | \leq s} \left \| \frac{\partial^{\alpha_1}\tilde{u}(r) }{\partial r^{\alpha_1}}\frac{n^{\alpha_2} }{r^{\alpha_2}}   \right\|_{L^2_\gamma}. \]
This definition comes from viewing 
$u(r, \theta) = \tilde{u}(r) \rme^{\rmi n \theta} \in h^s_{n,\gamma}$ 
as an element in $H^s_\gamma(\R^2)$, 
and writing the norm of this space using polar coordinates
\[ \|u\|_{H^s_\gamma} = \sum_{| \alpha | \leq s} \left \| \frac{\partial^{|\alpha|}u }{\partial r^{\alpha_1} \partial \theta^{\alpha_2}}\frac{1}{r^{\alpha_2}}   \right\|_{L^2_\gamma}. \]
Notice as well that Lemma \ref{l:polar2} then allow us to define an equivalent norm in $H^s_\gamma(\R^2)$, namely
\[ \|u\|^2_{H^s_\gamma} = \sum_n  \| u\|^2_{h^s_{\gamma,n}} .\]
We will use this definition in Subsection \ref{ss:fredholmProp} to prove the invertibility for some convolution operators.
\end{remark}


\subsection{Fourier Transform}
Here we recall some results from \cite{stein1971} regarding 
the direct sum decomposition of $L^2(\R^2)$ presented above, 
i.e. $L^2(\R^2) = \oplus h_{0,n}^0$. 
In particular, the next lemma shows that the spaces $h^0_{0,n}$ 
are invariant under the Fourier Transform, $\mathcal{F}$. 
To simplify notation, from now on we let $h_n = h^0_{0,n}$.
We also use the notation $L_r^2(\R^2)$ to describe the set of 
radially symmetric functions in $L^2(\R^2)$.

\begin{lemma}\label{l:invarianceF}
The Fourier Transform maps the spaces  
$$h_n = \{ f \in L^2(\R^2) \mid f( z) = g(r) \rme^{\rmi n \theta} , g \in L^2_r(\R^2)\}$$ 
back to themselves. In particular, given 
$f( z)= f(r \rme^{\rmi \theta}) = g(r) \rme^{\rmi n \theta} \in h_n$, 
then the Fourier transform of these functions can be written as
\[ \mathcal{F}[f(z)] = \mathcal{P}_n[g](\rho) \rme^{\rmi n \phi} = \breve{g}(\rho) \rme^{\rmi n \phi},\]
where
\[ \mathcal{P}_n[g] (\rho) =(-\rmi)^{n} \int_0^\infty g(r) J_n(r \rho) r \;dr, \]
and $J_n(z)$ is the $n$-th order Bessel function of the first kind.
Moreover,
\[ \mathcal{F}^{-1}[\hat{f}(w)] =  \mathcal{P}^{-1}_n[\breve{g}](r) \rme^{\rmi n \theta} = g(r) \rme^{\rmi n \theta},\]
with
\[ \mathcal{P}^{-1}_n[\breve{g}] (r) = \rmi^{n} \int_0^\infty \breve{g}(\rho) J_n(r \rho) \rho \;d\rho. \]

\end{lemma} 

The results of this lemma follow from the fact that the Fourier transform
 commutes with orthogonal transformations. 
 A detailed proof can be found in  \cite[Theorem 1.6]{stein1971},
 but we also provide a summary in the Appendix.
 
\section{The Convolution Operator}\label{s:convolution}
In this section we recall the assumptions made on the convolution kernels, $K$,
and summarize some of their properties.

\begin{hypothesis}\label{h:analytic}
The convolution kernel $K$ has a radially symmetric Fourier symbol 
$\hat{K}(\xi) = \hat{K}(|\xi|)$. As a function of $\rho = |\xi|$, the symbol $\hat{K}(\rho)$ 
can be extended to a uniformly bounded and analytic function on a strip 
$\Omega = \R \times ( -\xi_0, \xi_0) \subset \mathbb{C} $, for some constant $\xi_0>0$.
\end{hypothesis}

\begin{hypothesis}\label{h:taylor}
The symbol $\hat{K}(\rho)$ is symmetric and has a simple zero, 
which we assume is located at the origin $\rho=0$. 
This zero is of order $\ell =2$ and thus $\hat{K}(\rho)$  
has the following Taylor expansion near the origin:
\[ \hat{K}(\rho) \sim - \alpha \rho^{2} + \rmO( \rho^{4}) \quad \alpha>0.\]
\end{hypothesis}

The first result of this section is Theorem \ref{t:fredholmK}, 
which was proved in \cite{jaramillo2018, jaramillo2019} and
shows that these convolution operators are Fredholm when viewed as 
operators between Kondratiev spaces. 
Then, in Subsection \ref{s:diagonal} we prove 
 that the convolution with a radially symmetric function 
maps the spaces $h_n$ back to themselves, (see Lemma \ref{l:diagonal}).
This last result then implies that the operator $K\ast$ 
is a diagonal operator when we view its domain as a subset of 
$L^2_{\gamma}(\R^2) = \oplus h_n$.


\subsection{Fredholm properties}
To understand the need for Kondratiev spaces in establishing the Fredholm properties of $K$, 
consider first the pseudodifferential operator $( \mathrm{Id}- \Delta)^{-1} \Delta$ as a map from its domain 
$D \subset L^2(\R^2)$  back to $L^2(\R^2)$. 
This operator is the composition of the invertible map $(\mathrm{Id}-\Delta)^{-1}$, 
and the Laplacian. 
It has a zero eigenvalue embedded in its essential spectrum, 
and as a result one can use one of the corresponding eigenfunctions 
to construct Weyl sequences. 
These sequences then show that the map does not have closed range 
and therefore it is not a Fredholm operator.

To see this more clearly, consider for example just the Laplace operator, 
$\Delta :H^2(\R^2) \longrightarrow L^2(\R^2)$. 
Its kernel is spanned by harmonic polynomials, and although none of these functions are in $H^2(\R^2)$, 
one can use them to construct the following sequence: 
take $u_n = \chi( |x|/n) p(x,y)$, where $p(x,y)$ represents a harmonic polynomial and 
$\chi(|x|)$ is a smooth radial function equal to one when  $|x| <1$, 
and  equal to zero when $|x|>2$. 
Notice that this sequence does not converge in $H^2(\R^2)$. 
However $\| \Delta u_n \|_{L^2} \rightarrow 0$ as $n \rightarrow \infty$,  
showing that the operator does not have a closed range.

On the other hand, if we consider 
$\Delta : M^{2,2}_{\gamma-2}(\R^2) \longrightarrow L^2_{\gamma}(\R^2)$ 
and set $\gamma$ to be a large positive number, 
the above sequence would not be a Weyl sequence. 
Indeed, the algebraic decay imposed by the weight means that  
$\| \Delta u\|_{L^2_{\gamma}} \not \rightarrow 0$. 
In contrast, if we impose algebraic growth by picking $\gamma<1$, 
the above sequence would now converge to an element in the domain 
$M^{2,2}_{\gamma-2}(\R^2)$.

This heuristic argument justifies the results of the next theorem.

 \begin{theorem}\label{t:fredholmK}
Let  $\gamma \in \R $ and suppose the convolution operator 
$K: M^{2,2}_{\gamma-2}(\R^2) \longrightarrow H^2_\gamma(\R^2)$ 
satisfies Hypothesis \ref{h:analytic} and Hypothesis \ref{h:taylor}. Then, 
\begin{itemize}
\item if $1 + m < \gamma < 2 + m $ with $m \in \N$, the operator is Fredholm, injective, and has cokernel 
\[ \cup_{j=0}^m \mathcal{H}_j\]
\item if $ - m  < \gamma < 1- m $ with $m \in \N$, the operator is Fredholm, surjective, and has kernel 
\[ \cup_{j=0}^m \mathcal{H}_j\]
\end{itemize}
where $\mathcal{H}_j$ denotes the set of harmonic polynomials of degree $j$. 
On the other hand, if $ \gamma = m$ for some $m \in \N$, 
then the convolution operator does not have closed range.
\end{theorem}

The above result follow from Lemmas \ref{l:decomposition}, 
Lemma \ref{l:Minvertible} and Proposition \ref{p:fredholm}, 
which show that these convolution operators can be written as the 
composition of an invertible operator and a Fredholm operator. 
For the 1-d case the proof of Lemmas  \ref{l:decomposition},
and  \ref{l:Minvertible}
appear in  \cite{jaramillo2019} as Lemmas 3.9 and 3.10, respectively; while 
[Proposition 3.7, \cite{jaramillo2019}] is the analogue of our Proposition \ref{p:fredholm}
appearing below.
For the 2-d case, Lemmas 4.2 and 4.3
 in \cite{jaramillo2018} correspond to Lemmas  \ref{l:decomposition},
and  \ref{l:Minvertible} stated below, whereas Proposition  \ref{p:fredholm} in this manuscript appears as
 Proposition 3.6 in \cite{jaramillo2018}.

\begin{lemma}\label{l:decomposition}
The Fourier symbol $\hat{K}$ satisfying Hypotheses \ref{h:analytic} and \ref{h:taylor} 
admits the following decomposition
\[ \hat{K}(\xi ) =  \hat{M}_L(|\xi|)\hat{L}_{NF}(\xi) = \hat{L}_{NF}(\xi)  \hat{M}_R(|\xi|),\quad \xi \in \C^2, \]
where $\hat{L}_{NF}(\xi) = -|\xi|^2/(1 + |\xi|^2)$. Moreover, the symbols 
$\hat{M}_L(|\xi|), \hat{M}_R(|\xi|)$  together with their inverses are 
analytic and uniformly bounded functions of 
$\rho = |\xi|$, for $\rho \in \Omega\subset \C$ 
(see Hypothesis \ref{h:analytic} for the definition of $\Omega$).
\end{lemma}

Notice that because the Fourier symbols $\hat{M}_L, \hat{M}_R$, their inverses, 
and all their derivatives are analytic and uniformly bounded, 
it follows from Plancherel's Theorem that the corresponding operators 
$\mathcal{M}_{L/R} : H^s_\gamma(\R^2) \longrightarrow H^s_\gamma(\R^2)$, 
with $s \in \N \cup \{0\}$ and defined by
\[
\begin{array}{c c c}
H^s_\gamma(\R^2) & \longrightarrow &H^s_\gamma(\R^2)\\
u & \longmapsto & \mathcal{F}^{-1}( \hat{M}_{L/R} \hat{u}),
\end{array}
\]
are isomorphisms if $\gamma \in \Z_+$. 
This result can then be extended to values $\gamma \in \Z_-$ using duality, 
and to general $\gamma \in \R$ via interpolation, giving us the following lemma.

\begin{lemma}\label{l:Minvertible}
Given $s \in \N \cup \{0\}$, the operator 
$\mathcal{M}_{L/R}: H^s_\gamma(\R^2) \longrightarrow H^s_\gamma(\R^2)$, 
with Fourier symbol $\hat{M}_{L/R}(\xi)$  is an isomorphism for all $\gamma \in \R$.
\end{lemma}

The two lemmas above show us that the the convolution operators 
considered here are the composition of an invertible operator, 
 $\mathcal{M}_{L/R}$, and the pseudodifferential operator 
$ (\rmId - \Delta)^{-1}\Delta$. 
Therefore, the operators $K$ and $(\rmId - \Delta)^{-1}\Delta$
share the same Fredholm properties.

Now, to establish the Fredholm properties of the pseudodifferential operator
$$(\rmId - \Delta)^{-1}\Delta : M^{2,p}_{\gamma-2}(\R^2) \longrightarrow W^{2,p}_\gamma(\R^2),$$ 
one first notices that 
$(\rmId-\Delta): W^{s,p}_\gamma(\R^2) \longrightarrow W^{s-2,p}_\gamma(\R^2)$ 
can be written as a compact perturbation of 
$(\rmId-\Delta): W^{s,p}(\R^2) \longrightarrow W^{s-2,p}(\R^2)$, 
and is therefore invertible, see also \cite{jaramillo2019, jaramillo2018}. 
Then, in reference \cite{mcowen1979} it is shown that the the Laplacian 
$\Delta: M^{2,p}_{\gamma -2}(\R^2) \rightarrow L^p_\gamma(\R^2)$,
is  Fredholm. Combining these two results then leads to the next proposition.

\begin{proposition}\label{p:fredholm}
Let $p$ and $q$ be conjugate exponents, let $\gamma \in \R $, and consider the operator
\[ (Id-\Delta)^{-1}\Delta: M^{2,p}_{\gamma-2}(\R^2) \longrightarrow W^{2,p}_\gamma(\R^2),\]
then
\begin{itemize}
\item if $2/q + m < \gamma < 2/q + m +1$ with $m \in N$, the operator is Fredholm, injective, and has cokernel 
\[ \cup_{j=0}^m \mathcal{H}_j\]
\item if $2- 2/p - m -1 < \gamma < 2- 2/p - m $ with $m \in N$, the operator is Fredholm, surjective, and has kernel 
\[ \cup_{j=0}^m \mathcal{H}_j\]
\end{itemize}
where $\mathcal{H}_j$ denotes the set of harmonic polynomials of degree $j$. 
On the other hand, if $ \gamma = m$ for some $m \in \N$, then 
$\Delta$ does not have closed range.
\end{proposition}
The results of Theorem \ref{t:fredholmK} then immediately follow.

Finally, taking into account Lemma \ref{l:decomposition}, we make one further assumption on the convolution operator $K $.

\begin{hypothesis}\label{h:decomposition}
Given a Fourier symbol $\hat{K}$ satisfying the conditions in Lemma \ref{l:decomposition}, we have that 
\[ \hat{K}(\rho ) = \hat{M}_L(\rho) \frac{ -\rho^2}{1 + \rho^2} = \frac{ -\rho^2}{1 + \rho^2}\hat{M}_R(\rho), \]
where $\rho = | \xi|^2$. We assume that the analytic and uniformly bounded symbols, $\hat{M}_{L/R}$, satisfy  $$\hat{M}_L(\rho) = \hat{M}_R(\rho) = \hat{M}(\rho) = c( 1+ \hat{G}(\rho)),$$ with 
$c \in \R$, and $\hat{G}(\rho) \sim \rmO(1/\rho)$ as $\rho \to \infty$.
\end{hypothesis}


\subsection{Diagonalization}\label{s:diagonal}
In this subsection we prove that the convolution operators 
considered here map the spaces $h_n$ back to $h_n$. 
More precisely, given 
$u(z) = u_n(r) \rme^{\rmi n\theta} \in \mathcal{S} \cap h_{n}$, 
where $\mathcal{S}$ denotes the Schwartz space of rapidly decaying functions, 
we have that
\[ K \ast u = f(r) \rme^{\rmi n \theta} \quad \mbox{with} \quad f(r) = K_n \; \tilde{\ast} \; u_n.\]
In particular, $K_n$ is an appropriate radial function and the symbol 
$\tilde{\ast}$ denotes a convolution type of operator.
In other words, a Fourier series expansion in the angular variable, 
denoted here by $\mathcal{FS}$, diagonalizes the operator:
\[ \mathcal{FS}[K \ast u ] = \sum_n(K_n \tilde{\ast} u_n) \rme^{\rmi n \theta}. \]
In the next Subsection we will use this result to infer Fredholm properties 
for the restriction of the convolution operator $K$ to the subspace $h_n$.
 This will then allow us to pick a critical mode $n_0$ and the corresponding
  subspace where the normal form can be constructed. 

\begin{lemma}\label{l:diagonal}
Let $K$ be a radially symmetric kernel.
Then, the convolution with this kernel leaves the subspaces 
$\mathcal{S} \cap h_n=\{ u \in L^2(\R^2) \cap \mathcal{S} \mid u(r \rme^{\rmi \theta} )= \bar{u}(r) \rme^{\rmi n \theta} \quad \bar{u} \in L^2_r(\R^2)\} $ invariant.
\end{lemma}

\begin{proof}
First, notice that since $u \in \mathcal{S}$ the expression $K \ast u$ is well defined.
Second, in Lemma \ref{l:invarianceF} we proved that the Fourier transform
 leaves the spaces $h_n$ invariant. As a result the following diagram, 
 where $\mathcal{F}$ represents the Fourier Transform and $\mathcal{FS}$ 
 represents the Fourier series expansion on the angular variable, commutes.
\begin{center}
\begin{tikzpicture}
\draw [thick, ->] (0,2.5) -- (3,2.5);
\node[left] at (0,2.5)  {$L^2(\R^2)$};
\node[right] at (3,2.5)  {$L^2(\R^2)$};
\draw [thick, ->] (-0.5,2) -- (-0.5,0);
\draw [thick, ->] (3.5,2) -- (3.5,0);
\node[left] at (0,-0.25) {$\oplus h_n$};
\node[right] at (3,-0.25) {$\oplus h_n$};
\draw [thick, ->] (0,-0.25) --( 3,-0.25);
\node[above] at (1.5, -0.25) {$\mathcal{F}$};
\node[above] at (1.5, 2.5) {$\mathcal{F}$};
\node[left] at (-0.5, 1) {$\mathcal{FS}$};
\node[right] at (3.5, 1) {$\mathcal{FS}$};
\end{tikzpicture}
\end{center}

 The result now follows from our assumption that the kernel $K$ 
 is a radial function and therefore it has a radially symmetric Fourier symbol.  
 Indeed, we can see that 
 \[ \mathcal{FS}\;\Big [ \mathcal{F} [ K \ast u ] \Big ]  = \mathcal{FS}\;[ \hat{K}(|\xi|) \hat{u}(\xi)]  = \sum_n \left[ \hat{K}(|\xi|) \hat{u}(\xi) \right]_n \rme^{\rmi n \phi}  =  \sum_n \hat{K}(\rho) \hat{u}_n(\rho)  \rme^{\rmi n \phi},   \]
 where $\xi = \rho \rme^{\rmi \phi}$ and
 \[\hat{u}_n(\rho) = \frac{1}{2\pi} \int_0^{2\pi} \hat{u}(\xi) \rme^{-\rmi n \phi} \;d\phi.\]
The diagram then implies that,
\[ \mathcal{F}^{-1} [  \sum_n \hat{K}(\rho) \hat{u}_n(\rho) \rme^{\rmi n \phi} ] = \sum_n \mathcal{P}^{-1}_n[\hat{K}(\rho) \hat{u}_n(\rho)]  \rme^{\rmi n \theta}, \]
where $\mathcal{P}^{-1}_n$ is defined in Lemma \ref{l:invarianceF}. In other words,
$$ \mathcal{FS}[ K \ast u] = \sum_n (K_n \; \tilde{\ast} \; u_n)(r) \rme^{\rmi n \theta},$$
with
\[(K_n \; \tilde{\ast}\; u_n)(r) = \mathcal{P}^{-1}_n[\hat{K}(\rho) \hat{u}_n(\rho)].\]
\end{proof}

\subsection{ Fredholm Properties Revisited}\label{ss:fredholmProp}
In this subsection we summarize the Fredholm properties 
of the convolution operators $K \ast $ and $K \ast  + \rmi cn $ 
when considered as operators on the subspaces 
$ m^2_{\gamma,n}$ and $h^s_{\gamma,n}$, respectively.

\begin{lemma}\label{l:fredholmrevisited}
Let $\gamma \in \R $, $n \in \Z$, and consider the convolution kernel $K$ 
satisfying Hypotheses \ref{h:analytic} and \ref{h:taylor} 
restricted to the subspace 
$$m^2_{\gamma-2,n}= \{ u \in M^{2,2}_\gamma(\R^2) \mid u(r, \theta) = u_n(r)\rme^{\rmi n \theta}, \; \; u_n(r) \in M^{2,2}_{r,\gamma}(\R^2;\C), \; \bar{u}_n = u_{-n} \}.$$
Then,
$$K : m^2_{\gamma-2,n} \longrightarrow h^2_{\gamma,n}$$ 
is a Fredholm operator and
\begin{itemize} 
\item for $1-|n| < \gamma < |n|+1$, the map is invertible;
\item for $\gamma>|n|+1$ the map is injective with cokernel  spanned by $r^n \rme^{\rmi n \theta}$;
\item for $\gamma< 1-|n|$ the map is surjective with kernel spanned by $r^n\rme^{\rmi n \theta}$.
\end{itemize}
On the other hand, the operator is not Fredholm for integer values of $\gamma$.
\end{lemma}

\begin{proof}
First recall the results from Theorem \ref{t:fredholmK}, which show that the operator 
\[ K: M^{2,2}_{\gamma-2}(\R^2) \longrightarrow H^2_\gamma(\R^2)\]
is Fredholm for non integer values of $\gamma$.
Because the Fourier symbol for the convolution kernel $K$ is a radial function, 
Lemma \ref{l:diagonal} together with Theorem \ref{t:fredholmK} 
then show that the restriction operator
$$K: m^2_{n,\gamma-2} \longrightarrow h^2_{n,\gamma}$$
is not only well defined, but also Fredholm.

Finally, to obtain the description of the kernel and cokernel 
given in this Lemma, one can complexify 
$\R^2$, i.e. let $z = x + \rmi y$. 
Then the harmonic polynomials, 
which are the elements in the kernel and cokernel of $K$, 
are given by the real and imaginary parts of 
$(x + \rmi y)^n = z^n = r^n \rme^{\rmi n \theta}$.
\end{proof}

In the next section we will use the following 
Lemma which establishes the invertibility of the convolution operators 
$K  + \rmi cn $, restricted to the subspace 
$$h^s_{\gamma,n} =  \{ u \in H^s_\gamma(\R^2) \mid u(r, \theta) = 
u_n(r)\rme^{\rmi n \theta} \quad u_n(r) \in H^s_{r,\gamma}(\R^2;\C), \quad \bar{u}_n = u_{-n}\}.$$

\begin{lemma}\label{l:invertibleBn}
Let $s \in \N \cup \{0\}$, $\gamma \in \R$, $c \in \R\setminus\{0\}$ 
and consider the convolution kernel $K$ satisfying Hypotheses \ref{h:analytic}, \ref{h:taylor} and \ref{h:decomposition}.
Then, for all  $n \in \Z \setminus\{0\}$, the operator 
$L_n: h^s_{\gamma,n} \longrightarrow h^s_{\gamma,n}$
defined by
\[ L_n \;u(r)\rme^{\rmi n \theta}  = K \ast u(r) \rme^{\rmi n \theta} + \rmi cn \;u(r)\rme^{\rmi n \theta}\]
is invertible. Moreover,
\[  \| L_n u \|_{h^s_{\gamma,n}} \leq n C(\gamma) \|u\|_{h^s_{\gamma,n}} \qquad \mbox{and}\qquad \| L^{-1}_n f \|_{h^s_{\gamma,n}} \leq \frac{ \bar{C}(\gamma) }{n} \| f\|_{h^s_{\gamma,n}},\]
where $C(\gamma)$ and $\bar{C}(\gamma)$ are positive constants.
\end{lemma}

\begin{proof}
First consider the case of $\gamma \in \N \cup \{0\}$. 
Because the kernel is a radial function, 
Lemma \ref{l:diagonal} shows that the operator $K $ maps 
$h^s_{n,\gamma} \subset h_n$ to the space $h_n$. 
At the same time, using Remark \ref{r:normH} and Plancherel's theorem, 
given $u = u_n(r) \rme^{\rmi n \theta} \in h^s_{n,\gamma}$ 
with Fourier transform $\hat{u} = \breve{u}_n(\rho) \rme^{\rmi n \phi} \in h^\gamma_{n,s}$,
we have that for some generic constant $C>0$, 
\[ C  \| u_n \|_{h^s_{n,\gamma}} = \| u \|_{H^s_\gamma} = \| \hat{u} \|_{H^\gamma_s} = C \| \breve{u}_n \|_{h^\gamma_{n,s}}.\]

The results of the lemma then follow
if we show that the symbol $\hat{K}(\rho) + \rmi cn$, 
its inverse, and all their derivatives are uniformly bounded as functions of $\rho \in \R$.

From Hypothesis \ref{h:analytic} we know that as a function of $\rho = | \xi|$,
 the symbol $\hat{K}(\xi)= \hat{K}(|\xi|)$ is analytic on a strip $\Omega \subset \C$,
so that there exists a subdomain $\bar{\Omega}\subset \Omega \subset \C$ where 
$\hat{L}_n = \hat{K} (\rho) + \rmi cn$ is also analytic. 
Lemma \ref{l:decomposition} and Hypothesis \ref{h:decomposition}, then imply that this same symbol satisfies
$$\hat{K}(\rho) = \dfrac{ - \hat{M}(\rho) \rho^2}{1 + \rho^2},$$  
where $\hat{M}(\rho)$ is an analytic function that, 
together with its inverse and all its derivatives, is uniformly bounded on $\Omega \subset \C$.
Therefore, if we restrict $\rho \in \R$, then
\begin{align*}
 |\hat{L}_n(\rho)|^2 \leq & \sup_{\rho \in \R} \left(\frac{ - \hat{M}(\rho) \rho^2}{1 + \rho^2}\right)^2 + (cn)^2 < C + (cn)^2,\\
 |\hat{L}_n(\rho) |^2 \geq & \inf_{\rho\in \R} \left(\frac{ - \hat{M}(\rho) \rho^2}{1 + \rho^2}\right)^2 + (cn)^2 >  (cn)^2,
 \end{align*}
for some constant $C$. As a result, we also find that  $ |\hat{L}^{-1}_n(\rho)| < 1/|cn|$.

Straightforward calculations also show that all derivatives 
$D^\alpha \hat{L}_n(\rho)$ and $D^\alpha \hat{L}^{-1}_n(\rho)$,
with $\alpha$ satisfying $\alpha \leq \gamma $,
 are uniformly bounded. In particular,
  $$|D^\alpha \hat{L}_n(\rho)|<C(\gamma)\quad \mbox{and} \quad |D^\alpha \hat{L}^{-1}_n(\rho)|< C(\gamma)/(cn)^{|\alpha|+1},$$
where again $C(\gamma)$ represents a generic constant that depends on $\gamma$.

This proves the results of this lemma for the case of positive integer values of $\gamma$.
One can then extend the results to non-integer values of $\gamma$ by interpolation, 
and to negative values of $\gamma$ by duality.
\end{proof}

The following theorem from Reed and Simon's book \cite{reed1975},  
together with Hypothesis \ref{h:decomposition}, allows us to show that the operator
$K \ast + \zeta: W^{2, \infty}(\R^2) \rightarrow W^{2,\infty}(\R^2)$
is invertible for any $\zeta \in \C$, with $\mathrm{Im}(\zeta) \neq 0$. 
This is proved in Lemma \ref{l:boundinfinity}.

\begin{theorem}[Reed and Simon, Theorem IX.13]\label{t:reed}
Let $f$ be in $L^2(\R^n)$. Then $\rme^{b|x|} f \in L^2(\R^n)$ for all $b<a$ if and only if
$\hat{f}$ has an analytic continuation to the set $\{ \xi : | \mathrm{Im} \xi| < a\}$ with the property that 
for each $\eta \in \R^n$ with $|\eta | <a$, $\hat{f}(\cdot + \rmi \eta) \in L^2(\R^n)$ and for $b<a$
\[ \sup_{|\eta| \leq b} \| \hat{f}(\cdot + \rmi \eta)\|_{L^2} < \infty.\]
\end{theorem}
Notice that if $\rme^{b|x|} f \in L^2(\R^n)$, then by H\"older's inequality we have that $f$ is in $L^1(\R^n)$.

\begin{lemma}\label{l:boundinfinity}
Let $s \in \N$, fix $\zeta \in \C$ with $\mathrm{Im}(\zeta) \neq 0$, and consider the convolution kernel $K$ satisfying Hypotheses \ref{h:analytic}, \ref{h:taylor},  and \ref{h:decomposition}. Then,
the operator $$K \ast + \zeta: W^{s,\infty}(\R^2) \rightarrow W^{s,\infty}(\R^2)$$ is invertible.
\end{lemma}

\begin{proof}
Hypothesis \ref{h:decomposition} allows us to write the Fourier symbol for $K$ as 
\[ \hat{K}(\rho) = c \left[ -1 + \frac{1 - \rho^2 \hat{G}(\rho)}{1 +\rho^2}\right] = c \left [ -1 + \hat{G}_1(\rho) \right],\qquad \rho = |\xi|, \quad \xi \in \R^2.\]

Because $\hat{G}(\rho) \sim \rmO(1/\rho)$ as $\rho \to \infty$, it follows that the term $\hat{G}_1$ above,
 viewed as a function of $\xi \in \R^2$, is in $L^2(\R^2)$. Because this function is also smooth in $\rho$, there is a constant 
 $\tau >0$, such that $\hat{G}_1$ can be continued analytically to the set $\{ \xi \in \C^2: | \mathrm{Im} \xi| < \tau \}$. In addition 
we can pick $\tau$ small enough so that $\hat{G}_1$ satisfies the conditions of Theorem \ref{t:reed} above. 
As a result, $\mathcal{F}^{-1}[\hat{G}_1]$ is an $L^1(\R^2)$ function. 
Using Young's inequality we can then conclude that $K \ast + \zeta$  defines a bounded operator 
from $W^{s,\infty}$ to this same space.

Next we look at $(K\ast +\zeta) ^{-1}$ which thanks to Hypothesis \ref{h:decomposition} has Fourier symbol
\[ (\hat{K} \ast + \zeta)^{-1}= \frac{1 + \rho^2}{ (\zeta - \hat{M}(\rho) ) \rho^2 + \zeta}. \]
To show that this symbol defines a bounded operator we rewrite it as
\[  (\hat{K} \ast + \zeta)^{-1} = \frac{1}{ (\zeta - \hat{M}(\rho) ) \rho^2 + \zeta} \; +\;  \frac{1}{\zeta-c} \;+ \; \frac{ \zeta \hat{Q}(\rho) /(\zeta -c)}{  (\zeta - \hat{M}(\rho) ) \rho^2 + \zeta},\]
where $\hat{Q}(\rho) = \left[ ( \zeta-c) \rho^2) - ( \zeta - \hat{M}(\rho)) \rho^2 - \zeta \right] / \zeta$.

Notice that because $\hat{M}(\rho) = c( 1- \hat{G}(\rho) )$, we have that $\hat{Q}(\rho) \to -1$ as $\rho \to \infty$. 
In addition, because the expression $(\zeta - \hat{M}(\rho) ) \rho^2 + \zeta$ is smooth with respect to $\rho \in \R$, 
we can find a small number $\tilde{\tau}>0$ such that $[(\zeta - \hat{M}(\rho) ) \rho^2 + \zeta]^{-1}$ is analytic on the strip
$\{ \xi \in \C^2: | \mathrm{Im} \xi| < \tilde{\tau} \}$ and satisfies the assumptions of Theorem \ref{t:reed}. 
As before, using Young's inequality we obtain that $(\hat{K} \ast + \zeta)^{-1}:W^{s,\infty}(\R^2) \rightarrow W^{s, \infty}(\R^2)$ is bounded.
 \end{proof}

Armed with the results from Lemma \ref{l:fredholmrevisited}, Lemma
\ref{l:invertibleBn}, and  Lemma \ref{l:boundinfinity} we are now ready to derive our normal form.

\section{Normal Form}\label{s:normalform}

In this section we derive a normal form for showing the existence of rotating wave solutions,
$U(r, \theta) = U(r, \vartheta + ct) $,
to  oscillatory systems with nonlocal coupling.
These solutions satisfy the steady state equation,
\[ 0  = K \ast U - c \partial_\theta U + F(U;\lambda) \quad U \in \R^2, \quad x \in \R^2,\]
with a reaction term $F(U;\lambda)$ that satisfies the following assumptions.
 \begin{enumerate}
 \item $F$ depends only on the variable  $U$ and not its derivatives.
 \item $F(0,\lambda) =0$ for all $\lambda \in \R$, and
 \item  $D_UF(0,0) = A_0$ has a pair of complex eigenvalues $ \nu =\pm \rmi \omega$.
\end{enumerate}

 Our work in this section is split as follows. 
 We first establish the notation we will be using throughout this section. 
  Then, in Subsection \ref{ss:linear_operator} we look at the above equation 
  and its linearization about the trivial state, $U=0$. 
  We show that there exists an appropriate space $X$ and a projection 
  $P:X \longrightarrow X_\parallel$ which diagonalizes this operator, 
  splitting it into an invertible  and a bounded map.  
  In Subsection \ref{ss:ms} we expand on the ideas presented in the introduction
  and start the multiple-scale analysis. In this subsection we also 
use the projection $P$ to carry out the Lyapunov-Schmidt reduction and split the steady state equation
 into a reduced equation and a complementary subsystem.
In Subsection \ref{ss:ift} we show that the nonlinear terms are well defined
in the chosen space, $X$, and prove the existence of solutions to the complementary system
 using the implicit function theorem.
 Finally, in Subsection \ref{ss:normalform} we use the reduced equation
  together with the projection $P$ to derive our normal form.
 
   {\bf Notation:}
  As mentioned in the introduction, because the system is close
  to a Hopf bifurcation and the parameter $\lambda$ is close to its critical value of zero, 
we may assume that our solutions exhibit multiple-scales. In other words, letting $\ep$ denote a small parameter,
we may establish
   fast and slow variables, which we denote by $r, t$ and $R= \ep r, T = \ep^2 t$, respectively.
   In addition, because we are interested in rotating waves, the different time scales can be
written directly into the solution, leading to the preliminary ansatz
$$U(r,\theta) =U(r,\vartheta +c t) =  U(r, \vartheta + c^*t + \ep^2\mu t ),$$
 where we let $c = c^* + \ep^2 \mu$. The value of $\mu$ is left as a free parameter and 
the value of $c^*$ is chosen so that given any $0 \neq n_0 \in \Z$ we have that $ \rmi c^*n_0 = \rmi \omega$, 
the eigenvalue of the matrix $A_0 = D_UF(0;0)$.

In this section we also split the reaction term, $F(U;\lambda)$, 
into its linear, $A$, and nonlinear part, $\tilde{F}(U;\lambda)$.
Our assumptions on $F$ imply that the map $A$ depends on the parameter $\lambda$. 
When this parameter is near its critical value of $\lambda =0$, 
we may expand $A$ and its eigenvalues $\nu$
as follows,
\[ A = A_0 + \lambda A_1(\lambda), \]
\[ \nu = \nu_0 + \lambda \nu_1(\lambda) \in \C,\]
with $\nu_0 = \pm \rmi \omega$. 
Here we also let $W_1, W^*_1$, denote the right and left eigenvectors of the matrix $A$
corresponding to the eigenvalue $\nu = \rmi \omega + \rmO(\lambda)$, 
and we choose them so that their inner product satisfies $\langle W^*_1, W_1 \rangle =1$. 
This also leads to the relation
\[ A_1 W_1 =  \nu_1(\lambda) W_1.\]

With the above notation, we may rewrite the steady state equation as follows
\begin{equation}\label{e:zero}
0= \underbrace{K \ast U   + A_0U - c^* \partial_\theta U}_{LU} +  \underbrace{\left[ - \ep^2 \mu \partial_\theta U +  \lambda A_1(\lambda)U + \tilde{F}(U;\lambda) \right]}_{\mathcal{N}(U,\lambda,\mu)}.
\end{equation}

In the next subsection we concentrate on the operator $L$.

\subsection{The Linear Operator}\label{ss:linear_operator}

Our goal in this subsection is to determine a base space $X$ and a splitting, 
$ X_\parallel \oplus X_\perp$ such that, given 
$U = U_\parallel + U_\perp \in X$, the operator $L$ can be written as
\[ LU =  \begin{bmatrix} L_\parallel & 0 \\ 0 & L_\perp\end{bmatrix} \begin{bmatrix} U_\parallel \\ U_\perp \end{bmatrix},\]
with $L_\perp: D_\perp \subset X_\perp \rightarrow Y_\perp$ an invertible operator
 and $ L_\parallel: D_\parallel \subset X_\parallel \rightarrow Y_\parallel$ a bounded operator.
 
To motivate our choice of $X$, we recall again that the convolution operator $K \ast$ 
leaves the spaces $h_n$ invariant. 
We therefore start by looking at the restriction of the linear operator to these subspaces. 
That is, we consider
\[ LU_n \rme^{\rmi n \theta} = \left(K \ast \quad + A_0 - \rmi c^*n \mathrm{I}\right) U_n \rme^{\rmi n \theta}.\]
Notice that if $n_0$ satisfies $c^*n_0 = \omega$, then 
the matrix $B_{n_0} = ( A_0 - \rmi c^*n_0 \;\mathrm{I})$  has a nontrivial kernel. 
A short computation also shows that for all other integers, n,
the matrices $B_n = (A_0 - \rmi c^*n\; \mathrm{I}),$ 
have nonzero eigenvalues, $\nu_{1,2} = -\rmi (c^*n \pm \omega)$, 
which thanks to Lemma \ref{l:invertibleBn} implies that  the operators 
$$K \ast + B_n: h^s_{\gamma,n} \times h^s_{\gamma,n} \longrightarrow h^s_{ \gamma,n} \times h^s_{\gamma,n}$$
are invertible.  This suggest that we consider the following splitting
\begin{eqnarray}\label{e:ansatz}
 U =& \overbrace{ W_{1, n_0} w_1(r)\; \rme^{\rmi n_0 \theta} + \overline{W}_{1,n_0}\overline{w}_1(r) \;\rme^{- \rmi n_0 \theta}}^{U_\parallel} \\ \nonumber
 &+ \overbrace{ W_{2,n_0} w_2(r)\; \rme^{\rmi n_0 \theta} + \overline{W}_{2,n_0}\overline{w}_2(r)\; \rme^{- \rmi n_0 \theta} + \sum_{n \neq \pm n_0} U_n(r) \rme^{\rmi n \theta} }^{U_\perp},
 \end{eqnarray}
where  $W_{1,n_0},W_{2,n_0}$ are the right eigenvectors of $B_{n_0}$ 
corresponding to eigenvalues $\nu =0$ and $\nu = -2\rmi \omega$, respectively. 
Notice that $W_{1,n_0}$ is the same as $W_1$, the eigenvector associated with the matrix $A_0$.

 With this information we can define the projection 
 $P: X \rightarrow X_\parallel$ given by
\begin{equation}\label{e:projection}
 PU = \frac{1}{2\pi} \int_0^{2 \pi} \langle W^*_{1,n_0}, U \rangle W_{1,n_0} \rme^{-\rmi n_0 \theta} \;d\theta + 
 \frac{1}{2\pi} \int_0^{2 \pi} \langle \overline{W}^*_{1,n_0}, U \rangle \overline{W}_{1,n_0}\rme^{\rmi n_0 \theta}\; d\theta,
 \end{equation}
where $ W^*_{1,n_0}, \overline{W}^*_{1,n_0}$ are the normalized left eigenvectors 
associated with the zero eigenvalue of the matrices $B_{n_0}$ and $B_{-n_0}$. 
Similarly we have the complementary projection $(\mathrm{I} - P): X \rightarrow X_\perp$. 

Now that we have a  vector decomposition for $U \in\R^2 $, 
we need to choose a Banach space, $X= X_\parallel \oplus X_\perp$, for these functions.
Notice that we could let 
$X_\parallel = m^2_{\sigma-2, n_0} \times m^2_{\sigma-2, -n_0}$, 
where each component is in the direction of 
$ W_{1,n_0}$ and $\overline{W}_{1,n_0}$, respectively, and take 
$X_\perp \subset \oplus h^2_{\gamma,n} \times \oplus h^2_{\gamma,n}$.  
These choices
would then allow us to show that, for an appropriate domain $D$, the linear operator 
$L:D \subset X \rightarrow \oplus h^2_{\gamma,n} \times \oplus  h^2_{\gamma,n}$ is Fredholm. 
Roughly speaking, this holds thanks to Lemma \ref{l:invertibleBn},
which implies that the operator $L: D  \subset X_\perp \longrightarrow X_\perp$ is invertible,
and thanks to Lemma  \ref{l:fredholmrevisited}, 
which proves that $L: X_\parallel \longrightarrow h^2_{\sigma,n_0} \times h^2_{\sigma,-n_0}$, is Fredholm.
We point out that this result is independent of the value of $\gamma$ and $\sigma$, 
so that the domain $D$ can be composed of either algebraically decaying 
or algebraically growing solutions, and still the operator $L$ would be Fredholm. 
However, in order to guarantee  that the nonlinearities are well defined in the space 
$\oplus h^2_{\gamma,n} \times  \oplus h^2_{\gamma,n}$, 
one then needs to impose algebraic decay on the elements of 
$D$ (i.e. pick $\gamma>-1, \sigma-2>-1$). 
As a result, solutions that are uniformly bounded would be excluded from this domain. 
In particular rotating spiral waves with non-vanishing amplitudes 
would not be part of $D$.
To get around this, in what follows we define spaces that allow us to decompose
our solutions into a uniformly bounded part and an algebraically decaying part. 



{\bf The Space $X_\perp$:}
We first define our base space $\mathcal{H}$ as the external direct sum 
of an algebraically weighted Sobolev space and the space of twice differentiable 
and essentially bounded functions, that is
$\mathcal{H} = H^2_\gamma(\R^2) \oplus W^{2,\infty}(\R^2).$ 
We  select values of $\gamma>-1$, 
which implies that
functions in $H^2_\gamma(\R^2)$ decay algebraically and thus 
 capture the near field behavior of solutions,
while
functions in  $ W^{2,\infty}(\R^2)$ encode their far field behavior.

Because the space $ W^{2,\infty}(\R^2) \subset  H^2_\delta(\R^2) $, for any fixed $\delta<-1$, 
we may use Lemma \ref{l:polar2} to decompose this space as a direct sum. 
That is $W^{2,\infty}(\R^2)  = \oplus b_n^2$, where
 \begin{align*}
 b_n^s= & \{ u \in  W^{s,\infty}(\R^2)  \mid u(r, \theta) = u_n(r) e^{ \rmi n \theta},\; \bar{u}_n = u_{-n}, \; \; \mbox{and} \; \; u_n(r) \in W^{s,\infty}(\R^2; \C) \},
    \end{align*}
Therefore, we may also write $\mathcal{H} = \oplus \mathfrak{H}_{\gamma,n}$, where $ \mathfrak{H}_{\gamma,n}= h^2_{\gamma,n} \oplus b_n^2.$
These considerations, together with
Plancherel's theorem and Remark \ref{r:normH}, 
then allow us to define a norm for the space $\mathcal{H}$, 
\[ \|w\|^2_{\mathcal{H}}  =\left [\sum    \| u_n\|^2_{h^2_{\gamma,n}} \right ]+\| v\|^2_{W^{2, \infty}},\]
which is valid for  all $w = u + v \in \mathcal{H}$.
Similarly, we consider the space
\[ \mathcal{H}_{n_0} = \oplus_{n \neq \pm n_0} \mathfrak{H}_{\gamma,n},\]
which is a closed subspace of $\mathcal{H}$, and therefore inherits the same norm, $\| \cdot \|_{\mathcal{H}}$.

Finally, to define the space $X_\perp$ we first change our coordinate system in the vector space $\R^2$, so that 
given $U = (U_1, U_2) \in \R^2$ we have that $U_1= U_1(x,y)$ is the component in the direction of the vector
$\{ W_{1,n_0}\}$ and $U_2 = U_2(x,y)$ is the component in the direction of $\{ W_{2,n_0}\}$. Then
\[ U \in X_\perp = \mathcal{H}_{n_0} \times \mathcal{H}.\]

As for the domain of $L_\perp$, we let $D_\perp = D_1 \times D_2 \subset \mathcal{H}_{n_0} \times \mathcal{H}$, where $D_i$ is the space of smooth functions closed under the norm
 \[ \|w \|^2_{D_i} = \left[ \sum ( 1+ n^2)  \| u_n\|^2_{h^2_{\gamma,n}} \right] +  \| v\|^2_{W^{2,\infty}} , \quad i \in \{1,2\}. \]

\begin{lemma}\label{l:invertibleL_perp}
  Consider the convolution operator, $K$,
   satisfying Hypotheses \ref{h:analytic}, \ref{h:taylor}, and \ref{h:decomposition}.
  Then, the operator $L_\perp: D_\perp \longrightarrow X_\perp$, defined as
  \[ L_\perp U = K \ast U  +A_0 U - c^* \partial_\theta U\]
 is invertible.
\end{lemma}
\begin{proof}

Because the Fourier symbol $\hat{K}$ is radially symmetric, 
by Lemma \ref{l:diagonal} the operator $L_\perp$ is a block diagonal operator 
when we view its domain as a subspace of $ \mathcal{H} \times \mathcal{H}$, or equivalently $ \oplus (\mathfrak{H}_{\gamma,n} \times \mathfrak{H}_{\gamma,n}) $.
Therefore,
we can focus on how the operator acts on 
$\mathfrak{H}_{\gamma,n} \times\mathfrak{H}_{\gamma, n}   $. 
In addition, one notices that for elements 
$U_n\rme^{\rmi n \theta} \in (D_\perp \cap \mathfrak{H}_{\gamma,n} \times\mathfrak{H}_{\gamma,n} )$, 
the operator takes the form 
\[ L_\perp U_n \rme^{\rmi n \theta}  = ( K \ast   +B_n) U_n \rme^{\rmi n \theta}, \]
where, thanks to the above change of coordinates in the definition of $X_\perp$, the matrices $B_n = A_0 - \rmi c^* n$ are diagonal with  eigenvalues 
$\nu = \rmi( c^*n \pm \omega)$, distinct from zero.
We can therefore simplify the analysis by looking at how the operator acts on elements 
$u \in \mathfrak{H}_{\gamma,n} = h^2_{\gamma,n} \oplus b_n^2$. 

Letting $u(r, \theta) = u_n(r) \rme^{\rmi n \theta} + v_n(r) \rme^{\rmi n \theta}$, 
with $u_n(r) \rme^{\rmi n \theta} \in h^2_{\gamma,n} $ and $v_n(r) \rme^{\rmi n \theta}  \in b_n^2$,
 and defining $L_n u= K \ast u  - \rmi( cn \pm \omega) u $, we may write     
\[ L_ \perp u   =  L_n u_n (r) \rme^{\rmi n \theta}  +L_n v_n(r) \rme^{\rmi n \theta}. \]
Lemma \ref{l:invertibleBn} then shows that the operator 
$L_n:  h^2_{ \gamma,n} \rightarrow h^2_{ \gamma,n}$ is invertible for all $\gamma \in \R$,
and that given $L_n u \rme^{\rmi n \theta} = f \rme^{\rmi n \theta}$, 
 the following bounds
\[ \| f \|_{h^2_{\gamma,n}}= \| L_n u \|_{h^2_{\gamma,n}} \leq n C(\gamma) \| u \|_{h^2_{\gamma,n}}\]
and
\[ \|u \|_{h^2_{\gamma,n}} =  \| L^{-1}_n f \|_{h^2_{\gamma,n}} \leq \frac{ \bar{C}(\gamma) }{|cn\pm \omega|} \| f \|_{h^2_{\gamma,n}},\]
hold for some generic constants $C(\gamma)$ and $\bar{C}(\gamma)$.
Therefore, 
\[ L_n: h^2_{\gamma,n}  \longrightarrow h^2_{\gamma,n}, \qquad \mbox{and} \qquad
L_n^{-1}: h^2_{\gamma,n}  \longrightarrow h^2_{\gamma,n} ,\]
are bounded operators. 

 Next, the results from Lemma \ref{l:boundinfinity} together with our definition of the spaces $X_\perp$ and  $\oplus b^2_n = W^{2,\infty}(\R^2)$ allow us to conclude that the operator $L_\perp$ is an isomorphism from  $\oplus b_n^2$ back to itself, and from $\oplus_{n \neq \pm n_0} b_n^2$ back to this same space.

Finally, if we now take $F = (F_1,F_2)  \in X_\perp$ and $U = (U_1,U_2)  \in D_\perp$ 
such that $L_\perp U = F$, we see that for $i \in \{1,2\}$.
\begin{align*}
 \| F_i\|^2_{\mathcal{H}} = &\left[ \sum  \| f^{(1)}_{n} \|^2_{h^2_{\gamma,n}} \right]  + \|f^{(2)}\|^2_{W^{2,\infty} }, \\
  \leq & C_1(\gamma) \left[ \sum  (1+ n^2) \| u_n \|^2_{h^2_{\gamma,n}} \right]  +C_2\| v \|^2_{W^{2,\infty}}, \\
  \leq &  C(\gamma) \|U_i\|^2_{D_i}, 
  \end{align*}
and  
\begin{align*}
 \| U_i\|^2_{D_i} =&\left[ \sum (1+n^2)  \| u_n \|^2_{h^2_{\gamma,n}} \right] +\| v \|^2_{W^{2,\infty}}, \\
\leq & \tilde{C}_1(\gamma) \left[ \sum \frac{(1+n^2)}{ |cn \pm w|^2}  \|f^{(1)}_n \|^2_{h^2_{\gamma,n}} \right] + \tilde{C}_2 \|f^{(2)}\|^2_{W^{2,\infty}}, \\
\leq &\tilde{C}(\gamma) \|F_i\|^2_{\mathcal{H}},
\end{align*}
as desired.
\end{proof}

{\bf The Space $X_\parallel$:} We now concentrate on the space $X_\parallel$, which we define as
\begin{equation}\label{e:parallel}
 X_\parallel = \mbox{span}\{ W_{1,n_0} u_{n_0}, \overline{W}_{1,n_0} u_{-n_0}\}.
\end{equation}
Here $W_{1,n_0}$ is given as in \eqref{e:ansatz} and the functions 
 $u_{\pm n_0} \in \mathfrak{H}_{\gamma,\pm n_0}$, with $\gamma >-1 $.
 Notice as well that $X_\parallel$ is a closed subspace of $\mathcal{H}$ and it inherits its norm.
We also define the range of $L_\parallel$ as $Y_\parallel = X_\parallel$.
The next lemma shows that the operator $L_\parallel: X_\parallel \longrightarrow X_\parallel$ is bounded. 

\begin{lemma}
 Consider the convolution operator, 
 $K$, satisfying Hypotheses \ref{h:analytic}, \ref{h:taylor}, and \ref{h:decomposition}. 
 Then, the operator $L_\parallel:X_\parallel \longrightarrow X_\parallel$, 
defined as
   \[ L_\parallel U = K \ast U  +A_0 U - c^* \partial_\theta U\]
is bounded.
\end{lemma}

\begin{proof}
By an appropriate change of coordinates we can take
$X_\parallel = \mathfrak{H}_{\gamma,n_0} \times \mathfrak{H}_{\gamma,-n_0}$.
Since $c^*n_0 = \omega$, we can then write our operator $L_\parallel$ as
  \[\begin{array}{c c c c }
 L_\parallel: &\mathfrak{H}_{\gamma,n_0} \times \mathfrak{H}_{\gamma,-n_0}     &\longrightarrow&\mathfrak{H}_{\gamma,n_0} \times \mathfrak{H}_{\gamma,-n_0}\\[2ex]
 &(u_+ ,u_-) & \longmapsto & (K \ast u_+,K \ast u_-)
 \end{array}
 \]
Thus, without loss of generality we can concentrate on functions 
$u \in \mathfrak{H}_{\gamma,n_0} = h^2_{\gamma,n_0} \oplus b^2_{n_0} $.
To simplify notation we also write $n$ instead of $n_0$ and consider functions  
$u = u_0 + v$, with 
$u_0 \in h^2_{\gamma,n}$ and $v \in b^2_n$. 

Hypotheses \ref{h:analytic} and \ref{h:taylor},
  together with Lemmas \ref{l:decomposition} and \ref{l:Minvertible}, 
  imply that the Fourier symbol $\hat{K}(\xi)$, along with all its derivatives, 
  are analytic and uniformly bounded functions. 
  Using Plancherel's Theorem, a similar analysis as that of Lemma \ref{l:invertibleBn} 
  then implies that $K\ast:h^2_{\gamma, n} \longrightarrow h^2_{\gamma, n}$ is bounded. 
  At the same time, a similar proof as in Lemma \ref{l:boundinfinity} then shows that   $K\ast:b^2_n \longrightarrow b^2_n$ is also bounded.
  
  This leads to
\begin{align*}
\|K \ast u \|_{\mathcal{H}} = & \|K \ast u_0 \|_{h^2_{\gamma,n}} + \| K \ast v\|_{W^{2,\infty}},\\
\leq &  C_1 \| u_0 \|_{h^2_{\gamma,n}} +  C_2 \| v \|_{W^{2,\infty}},\\
\leq & C \| u \|_{\mathcal{H}} .
\end{align*}
The result of the Lemma then follows directly.

\end{proof}

\subsection{Multiple-Scales}\label{ss:ms}
 In this subsection we continue with the multiple-scale analysis started in the introduction. 
 We assume our solutions depend on fast and slow variables that are independent of each other,
 derive a hierarchy of three equations, and use the projection $P$ defined in the previous subsection to split 
 the last equation into a reduced equation and a complementary system.
 
Assuming fast variables, $r$ and $t$, and slow variables, $R = \ep r$ and $T= \ep^2 t$, our preliminary ansatz
$U(r,\theta,R;\ep,\mu) = U(r, R; \vartheta + c^*t + \ep^2 \mu t)$ can be expanded as, 
\begin{equation}\label{e:expansion}
 U (r,\theta,R; \ep,\mu) = \ep U_1(\theta,R;\ep,\mu) + \ep^2 U_2(\theta, R;\ep,\mu)  + \ep^3 U_3(r,\theta;\ep,\mu),\end{equation}
with 
\begin{equation}\label{e:ansatzU1_2}
U_1(\theta,R;\ep,\mu)= W_1 w(R;\ep,\mu) \rme^{\rmi n_0 \theta} + \overline{W}_1 \overline{w}(R;\ep,\mu) \rme^{-\rmi n_0 \theta }.
\end{equation}
 Using the notation of subsection \ref{ss:linear_operator}, 
 we assume that 
 $$U_1 \in X_\parallel \subset \mathcal{H} \times \mathcal{H},$$
 and 
 $$U_{2,3} \in D_\perp \subset \mathcal{H}\times \mathcal{H},$$
 where $\mathcal{H} =  H^2_{\gamma}(\R^2) \oplus W^{2,\infty}(\R^2)$ and  $0< \gamma $.
 Remark that while $U_1,U_2$ depend on the slow coordinate $R$, 
 we take $U_3= U_3(r, \theta;\ep,\mu)$. This mimics the analysis done when using 
 center manifold theory to derive amplitude equations. 
 We are assuming that the term $U_3$ evolves faster in the spatial 
 direction than either $U_1$ or $U_2$. 
 
At this time we also determine how the scaling $R = \ep r$ 
affects the operation of convolution with the kernel $K$. 
Given that $L^2_\gamma (\R^2)  = \oplus h_{\gamma,n}$, 
 we may assume that
 $u(r,\theta) = u_n(\ep r) \rme^{\rmi n \theta}$, without loss of generality.
 Then, using Lemma \ref{l:invarianceF}, 
a straight forward calculation shows that the Fourier Transform of this function is
$\mathcal{F}[ u_n(\ep r) \rme^{\rmi n \theta} ]= \breve{u}_n(\rho/ \ep) \rme^{\rmi n \phi} / \ep^2$. 
Therefore, 
\[\begin{split}
( K \ast u)(r) & = \mathcal{F}^{-1}[ \hat{K}(\xi) \hat{u}(\xi) ] \\
&=  \mathcal{F}^{-1}[ \hat{K}(\rho) \breve{u}_n(\rho/ \ep) \rme^{\rmi n \phi} / \ep^2] \\
&  = \mathcal{P}^{-1}_n[ \hat{K}(\rho)\breve{u}_n(\rho/ \ep)]   \rme^{\rmi n \phi} / \ep^2.
\end{split} \]
More precisely,
 \begin{align*}
( K \ast u)(r) = &\mathcal{P}^{-1}_n[ \hat{K}(\rho)\breve{u}_n(\rho/ \ep)]   \rme^{\rmi n \phi} / \ep^2\\
=& \left[ \frac{ \rmi^n }{\ep^2}\int_0^\infty \hat{K}(\rho) \breve{u}_n(\rho/ \ep) J_n(r\rho) \rho \;d\rho \right] \; \rme^{\rmi n \phi} \\
=&  \left[ \rmi^n \int_0^\infty \hat{K}(\ep P) \breve{u}_n(P) J_n(\ep r P) P \;dP \right]  \rme^{\rmi n \phi}\\
= &  \ep^2  \left[ \rmi^n \int_0^\infty \hat{\tilde{K}}_\ep( P) \breve{u}_n(P) J_n(R P) P \;dP \right]  \rme^{\rmi n \phi}\\
= & \ep^2 (\tilde{K}_\ep \ast u) (R),
\end{align*} 
where we used the change of coordinates $P = \rho/\ep$ in the third line, 
and defined $\tilde{K}_\ep$ through its Fourier symbol $\hat{\tilde{K}}_\ep( P) =\frac{1}{\ep^2}\hat{K}( \ep P)$.

{\bf Taylor Expansion:} 
Here we look in more detailed at the nonlinearities $ \tilde{F}(U;\lambda)$.
If we Taylor expand these terms, we obtain
\[  \tilde{F}(U;\lambda) =  MUU + NUUU + \cdots,\]
where
 \begin{align*}  
  (MUU)_i = & \frac{1}{2!} \partial_{jk} \tilde{F}_i(0) U_jU_k,\\
  (NUUU)_i = & \frac{1}{3!} \partial_{jk\ell} \tilde{F}_i(0) U_jU_kU_\ell.
\end{align*}
To keep the nonlinearities as general as possible, 
we assume as well that each term in the series depends 
on the parameter $\lambda$ and has expansions of the form
\begin{align*}
 M(\lambda) = &M_0 + \lambda M_1(\lambda)  = M_0 + \ep^2 \bar{\lambda} M_1(\lambda), \\
 N(\lambda) = & N_0 + \lambda N_1(\lambda)  = N_0 + \ep^2 \bar{\lambda} N_1(\lambda),\quad \mbox{etc}.
 \end{align*}
{\bf Equating Coefficients:} 
For convenience, we again recall equation \eqref{e:zero},
\[
0= \underbrace{K \ast U - c^* \partial_\theta U  + A_0U}_{L} +  \underbrace{\left[ -  \ep^2 \mu \partial_\theta U +  \lambda A_1(\lambda)U + \tilde{F}(U;\lambda) \right]}_{\mathcal{N}(U,\lambda,\mu)}.
\]
Inserting the ansatz \eqref{e:expansion} 
into the above equation, letting
$\lambda = \ep^2 \bar{\lambda}, $
and collecting terms of equal order in $\ep$, 
gives us the next three relations.

 At $\rmO(\ep)$:
 \[ c^*\partial_\theta U_1 - A_0 U_1 =0.\]
 
 At $\rmO(\ep^2)$:
 \[ c^*\partial_\theta U_2 - A_0 U_2 =  M_0U_1U_1.\]
 
 And at higher orders:
\[\begin{split}
  c^*\partial_\theta U_3 - K \ast U_3 - A_0 U_3  = & - \mu( \partial_\theta U_1+  \ep \partial_\theta U_2 + \ep^2 \partial_\theta U_3) +   \tilde{K}_\ep \ast ( U_1 + \ep U_2 )  \\
  &+  \bar{\lambda} A_1(\lambda) [ U_1 + \ep U_2 + \ep^2 U_3] \\
  & +  \frac{1}{\ep^3} \left[ \tilde{F}(U;\lambda)- \ep^2 M_0 U_1U_1\right] .
\end{split}\]

We immediately notice that the first equation is satisfied if 
 $ U_1(R) = W_1 w(R) \rme^{\rmi n_0  \theta} + \overline{W}_1 \overline{w}(R) \rme^{-\rmi n_0 \theta},$
 where $W_1$ is the eigenvector for $A_0$ associated with the eigenvalue $\nu = \rmi \omega$.
This definition is consistent with our assumption that $U_1 \in X_\parallel$.
 Recall that this implies that in an appropriate coordinate system, 
 $w \rme^{\pm \rmi n_0 \theta} \in \mathfrak{H}_{\gamma, \pm n_0}  \subset H^2_\gamma(\R^2) \oplus W^{2,\infty}(\R^2)$
  with $0< \gamma $.

 To solve the second equation, notice that the right hand side involves the term 
 \[ U_1U_1 = W_1W_1 w(R)^2  \rme^{2 \rmi n_0 \theta} +2 W_1\overline{W}_1 |w|^2 + \overline{W}_1\overline{W}_1 \overline{w}(R)^2  \rme^{-2 \rmi n_0 \theta}.\]
 Thus, we conclude that $U_2$ must be of the form 
  \[ U_2 = V_1 w^2 \rme^{2\rmi n_0 \theta} + V_0 |w|^2 + V_{-1} \overline{w}^2 \rme^{-2\rmi n_0 \theta},\]
which leads to the next 3 linear equations for the vectors $V_1,V_0,V_{-1}$,
\begin{align*}
(2 \rmi n_0 c^* - A_0) V_1 = & M_0 W_1 W_1,\\
(-2 \rmi n_0 c^* - A_0) V_{-1} = & M_0 \overline{W}_1 \overline{W}_1,\\
-A_0 V_0 = & 2 M_0 W_1 \overline{W}_1.
\end{align*}

Because the function
$u = w \rme^{\pm \rmi n_0 \theta}  \in \mathfrak{H}_{\gamma,\pm n_0} $,
 Lemma \ref{l:nonlinearities} in the next subsection shows that terms of the form
$u^2,\bar{u}^2, |u|^2$, and in fact any power $u^p$, are in 
$\mathfrak{H}_{\gamma,\pm n_0}$. It then follows that $U_2$ is indeed in $X_\perp$.

 Finally, we use the projection 
 $P:X_\parallel \times D_\perp \longrightarrow X_\parallel$, 
 defined using \eqref{e:projection},  to split the third equation into the system
\begin{align}\label{e:reduced}
0 =& \tilde{K}_\ep \ast U_1 -  \mu \partial_\theta U_1 + \bar{\lambda} A_1(\lambda) U_1 + \frac{1}{\ep^3} P \left[ \tilde{F}(U;\lambda) - \ep^2 M_0U_1U_1 \right], \\[2ex] \label{e:ift}
0 =& - c^* \partial_\theta U_3 + K \ast U_3 + A_0 U_3 -  \mu( \ep \partial_\theta U_2 + \ep^2 \partial_\theta U_3) + \tilde{K}_\ep \ast (\ep U_2 + \ep^2 U_3) \\ \nonumber
&+ \bar{\lambda} A_1(\lambda) (\ep U_2 + \ep^2 U_3) + \frac{1}{\ep^3} ( \mathrm{I} - P) \left[ \tilde{F}(U;\lambda) - \ep^2 M_0U_1U_1 \right].
\end{align}

In the next subsection we show that the last equation defines an operator that satisfies
the conditions of the implicit function theorem. As a result solutions, $U_3$, to equation \eqref{e:ift} exist,
and they depend smoothly on $U_1, \ep,$ and $\mu$. 
In Subsection \ref{ss:normalform} we use this information, 
together with the reduced equation \eqref{e:reduced} and the projection $P$,
to derive our normal form.


\subsection{ Implicit Function Theorem:}\label{ss:ift}
We now look at the right hand side of equation  \eqref{e:ift} as an operator 
$$G_2(U_1, U_3; \ep,\mu): X_\parallel \times D_\perp \times \R^2 \longrightarrow X_\perp,$$ 
and prove that it satisfies the conditions of the implicit function theorem (recall that $U_2= U_2(U_1)$).
As a consequence, we obtain the existence of neighborhoods 
 $\mathcal{B} \subset \R^2$ and $\mathcal{U} \subset X_\parallel$, 
 with $(0,\mu^*) \in \mathcal{B}$ and $0 \in \mathcal{U}$, and a  map 
$\Psi : \mathcal{U} \times \mathcal{B} \longrightarrow  D_\perp $, 
such that $U_3 =  \Psi(U_1;\ep,\mu)$ satisfies 
\begin{align*}  
0 = &G_2(U_1, \Psi(U_1;\ep,\mu); \ep,\mu) \\
0 = & D_{U_1} \Psi(0;\ep,\mu)
   \end{align*}
   for all $U_1 \in \mathcal{U}$ and all $(\ep, \mu) \in \mathcal{B}$.

First, inspecting expression \eqref{e:ift} one can check that $G_2$ is smooth 
in all its variables and that given any $\mu=\mu^*$ it satisfies
 $G_2(0,0;0,\mu^*) =0$. 
In addition, the Fr\'echet derivative of $G_2$ evaluated at $U=0, \ep =0, \mu=\mu^*$ is given by
\[ D_UG_2(0,0;0,\mu^*) U =    K \ast U + A_0 U  - c^*\partial_\theta U, \]
 which is exactly the form of $L_\perp$ stated in Lemma \ref{l:invertibleL_perp}. 
 Therefore, $D_UG_2(0,0;0) : D_\perp \longrightarrow X_\perp$ 
 defines an isomorphism. We are left with showing that the operator is well defined.
 In particular, we need to show that the terms in the expression,
\[ \begin{split}
\mathcal{N} = & - \bar{\mu}( \ep \partial_\theta U_2 + \ep^2 \partial_\theta U_3) + \tilde{K}_\ep \ast (\ep U_2 + \ep^2 U_3) \\
&+ \bar{\lambda} A_1(\lambda) (\ep U_2 + \ep^2 U_3) + \frac{1}{\ep^3} ( \mathrm{I} - P) \left[ \tilde{F}(U;\lambda) - \ep^2 M_0U_1U_1 \right]
\end{split}\]
are in the space $X_\perp$.

First, because $U_2, U_3 \in D_\perp \subset X_\perp$ one can immediately see 
 that all linear terms
in the definition of $\mathcal{N}$ 
are well defined.
The results from Lemma \ref{l:nonlinearities}, 
which we state in the next paragraph, 
together with the projection $( \mathrm{I}-P)$, 
then show that all other higher order terms $\rmO((\lambda+ \mu)(U)^2)$ 
also map elements in 
$X_\parallel \oplus D_\perp \times \R^2$ to elements in $X_\perp$. 
Notice that Lemma \ref{l:nonlinearities} 
provides a more general result than what we need, 
and in the present argument we are using the fact that 
$D_\perp \subset X_\perp \subset \mathcal{H}\times \mathcal{H}$, as well as
$X_\parallel \subset \mathcal{H}\times \mathcal{H}$

\begin{lemma}\label{l:nonlinearities}
Let  $ \gamma \in \R$ with $0 < \gamma $ 
and let $p$ an integer such that $p\geq2$. Then, the map
\[
\begin{array}{c c c}
\tilde{N}:   \mathcal{H} \oplus \mathcal{H} & \longrightarrow & \mathcal{H}\\
(u_\parallel + u_\perp) & \longmapsto& (u_\parallel+ u_\perp)^p
\end{array}
\]
where $ \mathcal{H} =  H^2_\gamma(\R^2) \oplus W^{2,\infty}(\R^2)$,
 is well defined.
\end{lemma}
\begin{proof}
It is enough to show that $\mathcal{H}$
is  a Banach algebra. That is, given a function $w \in \mathcal{H}$, 
we need to show that any power $w^p $, with $p\geq 2$, belongs to this same space. 
If we let $w = u+v $, with $u \in H^2_\gamma(\R^2)$ and 
$v \in W^{2, \infty}(\R^2) $, we obtain the following expression for $w^p$,
\[ w^p = (u+ v)^p = v^p + \sum_{k=1}^p {p \choose k} v^{p-k} u^k.\]
Notice that $v^p \in W^{2, \infty}(\R^2)$, so that we are left with 
showing that the rest of the terms in the sum are in $H^2_\gamma(\R^2)$. 
Because $v^{p-k}$ is a bounded function for all $k \in [1,p] \cap \N$, 
we only need to show that $u^k$ is in $H^2_\gamma(\R^2)$ for all integers $k \geq 1$. 
To do this, we first prove that elements in $H^2_\gamma(\R^2)$ are uniformly bounded.

Since $u \in H^2_\gamma(\R^2)$, thanks to the Sobolev embeddings we have that
$|u(x)| \langle x \rangle^\gamma \in H^2(\R^2) \subset C_B(\R^2)$.
Then, because $\gamma>0$ we obtain that $ |u(x)| < \langle x \rangle^{-\gamma} <C$.
Therefore, $u(x)$ is a uniformly bounded function,
and it then follows that 
$u^k \in L^2_\gamma(\R^2)$ for any $k \geq 1$. 

Similarly, we find that the derivatives 
$D(u^k) = k u^{k-1} Du$ 
are well defined, since they are the product of a bounded function, $u^{k-1}$, with the 
$L^2_\gamma(\R^2)$ function $Du$. 
As for the second derivatives, 
$D^2(u^k) = k(k-1) u^{k-2} ( Du)^2 + k u^{k-1} D^2u$,  
this same argument shows that the last term is well defined, 
since it involves  the product of an $L^2_\gamma(\R^2)$ function, 
$D^2u$, with the bounded function $u^{k-1}$ .

 We are left with showing that the expression $( Du)^2$ is in 
$L^2_\gamma(\R^2)$. 
Here we can use  the Sobolev embedding,
$|Du(x)| \langle x \rangle^\gamma \in H^1(\R^2) \subset L^q(\R^2)$ for $ 2 \leq q < \infty$, 
together with H\"older's inequality to conclude that  
$\| (Du)^2 \|_{L^2_\gamma}$ is bounded. 
Indeed,
\[\begin{split}
\| (Du)^2 \|^2_{L^2_\gamma} = & \int_{\R^2} |Du|^4 \langle x \rangle^{2\gamma} \;dx \\
\leq & \left [  \int_{\R^2} \left( |Du| \langle x \rangle^{\gamma}\right)^2 \;dx  \right]^{1/2}\left [  \int_{\R^2} \left( |Du|^3 \langle x \rangle^{\gamma}\right)^2 \;dx  \right]^{1/2}\\
\leq & \left [  \int_{\R^2} \left( |Du| \langle x \rangle^{\gamma}\right)^2 \;dx  \right]^{1/2}\left [  \int_{\R^2} \left( |Du| \langle x \rangle^{\gamma/3}\right)^6 \;dx  \right]^{1/2}\\
\leq & \| Du \|_{L^2_\gamma(\R^2)} \| Du\|^3_{L^6_\gamma(\R^2)},
\end{split}\]
where the last inequality holds provided
$\langle x \rangle^{\gamma/3} < \langle x \rangle^{\gamma}$, i.e. $\gamma>0$.
This completes the proof.
\end{proof}

\subsection{Normal Form:}\label{ss:normalform}
Equation \eqref{e:reduced} will give us  our normal form. 
To simplify this expression we determine nonlinear terms up to 
order $\rmO(\ep^3)$ explicitly. Recall that 
$$\tilde{F}(U;\lambda) = MUU + NUUU + \cdots.$$
Using the notation from the start of this section, we find that
\begin{align*}
MUU = & (M_0 + \ep^2\bar{\lambda}M_1(\lambda))( \ep U_1 + \ep^2 U_2 + \ep^3 U_3)^2\\
MUU  = & \ep^2M_0 U_1U_1 + 2 \ep^3 M_0 U_1U_2 + \rmO(\ep^4), \\[3ex]
NUUU = & (N_0 + \ep^2\bar{\lambda}N_1(\lambda))( \ep U_1 + \ep^2 U_2 + \ep^3 U_3)^3\\
NUUU= & \ep^3 N_0 U_1U_1U_1 + \rmO(\ep^4).
 \end{align*}
 
 The reduced equation \eqref{e:reduced} is then given by
 \[\begin{split}
  \tilde{K}_\ep \ast U_1& -  \mu \partial_\theta U_1 + \bar{\lambda} A_1(\lambda) U_1+ \\
  &  P\left[ 2 M_0 U_1U_2+ N_0 U_1U_1U_1 + \rmO \Big (\ep (|U_1||U_3| + |U_2 + \ep U3|^4) \Big) \right]=0, 
  \end{split}\]
 which after projecting onto the space $X_\parallel$ results in the CGL-type equation
 \begin{equation}\label{e:cgl}
 0 =  \tilde{K}_{\ep,n_0} \ast w -  \mu \rmi n_0 w + \bar{\lambda} \nu_1(\lambda) w + (a_1+a_2) |w|^2 w + \rmO(\ep |w|^4w),
 \end{equation}
 and its complex conjugate. Here the symbol $\tilde{K}_{\ep,n_0}$ represents the action of the operator $ \tilde{K}_{\ep}$ on the space $X_\parallel $. The constants $a_1, a_2$ are found using the expressions for $U_1$ and $U_2$, via the relations,
 \begin{align*}
 a_1 = & \langle W^*_{1} , 2 M_0 (W_1V_0 + \overline{W}_1V_1) \rangle,\\
 a_2 = & \langle W^*_{1} , N_0 (W_1W_1 \overline{W}_1) \rangle.
 \end{align*}
We refer to equation \eqref{e:cgl}, including all higher order terms, as the normal form.


\section{Existence of Solutions and Validity of the Normal Form}\label{s:validity}
In this section we prove the existence of solutions to the normal form, i.e. equation  \eqref{e:cgl}. 
We use this result together with the analysis presented in Section \ref{s:normalform}
to show the validity of this equation.

Notice that our definition of the normal form includes all higher order terms, including those
summarized in the expression $\rmO(\ep|w|^4w)$.
Although we don't have a precise description for them,
we can still prove the existence of solutions to equation \eqref{e:cgl} 
using the implicit function theorem, provided these higher order terms are all well defined.
As shown in Proposition  \ref{p:existence}, which we recall and prove below, this is just consequence of Lemma \ref{l:nonlinearities}.
As a result we are able to establish
 the existence of solutions to this normal form, and show that they are valid in a 
neighborhood of $(\ep,\mu) = (0,\mu^*)$, where $\mu^*$ is an arbitrary number different from zero.

The above result, together with the Lyapunov-Schmidt reduction presented in 
Section \ref{s:normalform}, then imply the existence 
 of small amplitude solutions to the steady state equation \eqref{e:main_steady}.
These solutions take the form,
\[ U(r,\theta;\ep,\mu) = \ep U_1(\theta, \ep r; \ep, \mu) + \ep^2 U_2(\theta, \ep r; \ep, \mu)+ \ep^3 U_3(\theta, r; \ep, \mu).\]
Moreover, they are unique and valid in a small neighborhood of  $(\ep,\mu)=(0,\mu^*)$.
Consequently,  if  $w(R;\ep,\mu)$ is a solution to equation \eqref{e:cgl} and $U$ is a solution to 
the steady state equation, then the approximation
 \begin{equation}\label{e:ansatzU11}
  U_1(r, \vartheta; \ep,\mu) =  \ep (W_1 w(\ep r;\ep,\mu) \rme^{\rmi n_0 (\vartheta + (c^* +\mu) t )} + \overline{W}_1 \overline{w}(\ep r;\ep,\mu) \rme^{-\rmi n_0 (\vartheta + (c^* +\mu) t) })
  \end{equation}
 satisfies 
 $$\| U- U_1 \|_{C_B}<  \|\ep^2 U_2 + \ep^3 U_3 \|_{C_B}  < \ep^2.$$
  Here we used the fact that $U_{2,3} \in \mathcal{H} \subset C_B$ .
  Thus, Proposition \ref{p:existence} and the Lyapunov-Schmidt reduction 
  of Section \ref{s:normalform} also imply the validity of our  normal form equation. 
  Together, these results give us our main theorem, Theorem \ref{t:main}.
 \begin{customthm}{1}\label{t:main}
 Let  $\gamma>0$, $n_0 \in \Z$, and suppose $w(R;\ep,\mu) \in \mathfrak{H}_{\gamma,n_0} $ is a solution to equation \eqref{e:cgl}. 
Then, there exist unique solution $U(r, \theta) $ of the steady state equation \eqref{e:main_steady} and constants $ C,\ep_*>0$, such that for all $ \ep \in (-\ep,\ep_*)$ 
the estimate
\[ \| U(r,\theta) - U_1(r,\theta) \|_{C_B} < C\ep^2,\]
with $U_1$ as in \eqref{e:ansatzU11}, holds.
\end{customthm}
 
The rest of this section is dedicated to proving Proposition \ref{p:existence}.

\begin{customprop}{1.1}\label{p:existence}
Given real numbers $\mu^* \neq 0$,  $\gamma >0 $, 
and an integer $n$, there exists positive constants $\ep_0, \mu_0$, and a $C^1$ map 
\[
\begin{array}{ c c c c}
\Gamma: & (-\ep_0, \ep_0) \times ( \mu^*-\mu_0, \mu^*+ \mu_0) & \longrightarrow & \mathfrak{H}_{\gamma,n}\\
& (\ep,\mu) & \longmapsto & w(R; \ep, \mu)
\end{array}
\]
such that $w(R;\ep,\mu)$ is a solution to the equation
\begin{equation}\label{e:cgl2}
 0 =\tilde{K}_\ep \ast w + (\mu^* +\mu) \rmi n w + \lambda w + a |w|^2 w + \rmO(\ep |w|^4w).
\end{equation}
Here $\lambda \in \R, a \in \C$ are nonzero constants, 
and the Fourier symbol, $\ep^2 \hat{\tilde{K}}_\ep(\xi) = \hat{K}(\ep \xi)$ satisfies hypotheses \ref{h:analytic}, \ref{h:taylor},  and \ref{h:decomposition}.
\end{customprop}

To prove the proposition we first recall  $  \mathfrak{H}_{\gamma,n}=  h^2_{\gamma,n} \oplus b^2_n $,
and define the operators
\begin{align*}
 \mathcal{L}_{\mu,0} w &= \hat{M}(0) \Delta w + (\lambda +  (\mu^*+\mu) \rmi n ) w,\\
 \mathcal{L}_{\mu,\ep}w & = \tilde{K}_\ep \ast w + (\lambda +  (\mu^*+\mu) \rmi n ) w,
 \end{align*}
which we show below in Lemma \ref{l:Lmu} have bounded inverses with domain $\mathfrak{H}_{\gamma,n}$, 
and are also $C^1$ with 
respect to the parameters $\mu$ and $\ep$.
Notice that in the definition of $\mathcal{L}_{\mu,0}$ we have used the properties of
$\hat{\tilde{K}}_\ep(\xi) = \hat{K}(\ep \xi)/ \ep^2$, in particular Lemma \ref{l:decomposition},
to conclude that when $\ep =0$, the convolution with $\tilde{K}_\ep$ reduces to the Laplace operator.
The constant $\hat{M}(0)$, is just the Fourier symbol from Hypothesis \ref{h:decomposition} evaluated at zero.

Now, preconditioning  the normal form  with  $\mathcal{L}^{-1}_{\mu,\ep}$, 
we may view the right hand side of equation \eqref{e:cgl2} as an operator 
$F: \mathfrak{H}_{\gamma,n} \times \R \times \R  \longrightarrow \mathfrak{H}_{\gamma,n}$, 
given by
\[ F(w; \ep,\mu) =  \mathrm{I} w +  \mathcal{L}^{-1}_{\mu, \ep} \left[   a |w|^2 w + \rmO(\ep |w|^4w) \right].\]
The zeros of  $F$ then correspond to solutions of the equation, 
which we can find using the implicit function theorem. 

It is clear that the operator $F$ 
satisfies $F(0; 0,0) =0$, and that its Fr\'echet derivative 
$D_wF(0;0,0) = \mathrm{I}:\mathfrak{H}_{\gamma,n} \longrightarrow \mathfrak{H}_{\gamma,n} $  
defines an invertible operator. That $F$ is also well defined follows 
from Lemma \ref{l:nonlinearities} and Lemma \ref{l:Lmu}.  
Indeed, Lemma \ref{l:nonlinearities} shows that all nonlinearities of the form 
$|w|^p w^q$, for $q,p \in \N$, define a bounded map from 
$\mathfrak{H}_{\gamma,n}$ back to itself. 
On the other hand, Lemma \ref{l:Lmu}, which we state and prove next,
 shows that $\mathcal{L}_{\mu,\ep}^{-1}$ is bounded and
  that it is also continuously differentiable with respect to the parameters $\mu$ and $\ep$. 
As a result we also obtain that the operator $F$ is  continuously differentiable 
with respect to these parameters. We  may therefore apply the implicit function theorem, and the results of Proposition \ref{p:existence} then follow.

We now concentrate on proving the desired properties of the linear operator $\mathcal{L}_{\mu,\ep}$.

   \begin{lemma}\label{l:Lmu}
 Fix $\mu^*\neq 0, \lambda \neq 0 \in \R$, let $0 \neq n \in \Z,$  $\gamma \in \R$, and take $ \ep, \mu, $ to
 be real numbers. Consider as well the convolution kernel 
 $\tilde{K}_\ep$, with Fourier symbol $\ep^2 \hat{ \tilde{K}}_\ep(\xi) = \hat{K}(\ep \xi)$ satisfying Hypotheses \ref{h:analytic}, \ref{h:taylor}, and \ref{h:decomposition}, and define
 \begin{align*}
 \mathcal{L}_{\mu,0} w &= \hat{M}(0) \Delta w + (\lambda +  (\mu^*+\mu) \rmi n ) w,\\
 \mathcal{L}_{\mu,\ep}w & = \tilde{K}_\ep \ast w + (\lambda +  (\mu^*+\mu) \rmi n ) w,
 \end{align*}
Then, the inverse operators
\[
\begin{array}{c c c c }
\mathcal{L}_{\mu,0}^{-1} :& \mathfrak{H}_{\gamma,n} & \longrightarrow & h^4_{\gamma,n} \oplus b^2_n   \\[2ex]
\mathcal{L}_{\mu,\ep}^{-1} :& \mathfrak{H}_{\gamma,n}& \longrightarrow &\mathfrak{H}_{\gamma,n},\\[2ex]
\end{array}
\]
where $\mathfrak{H}_{\gamma,n} = h^2_{\gamma,n} \oplus b^2_n$,
are bounded.
 Moreover, 
$\mathcal{L}_{\mu,\ep}^{-1}$ is $C^1$ with respect to $\mu$ and $\ep$, for all $ \mu, \ep \in \R$.
 \end{lemma}

 \begin{proof}

 {\bf Step 1:} We first prove that  the operator $\mathcal{L}_{\mu,\ep}^{-1}$ is bounded.
 Because the Fourier symbol $\ep^2 \hat{ \tilde{K}}_\ep(\xi) = \hat{K}(\ep \xi)$ 
 satisfies Hypotheses \ref{h:analytic}, \ref{h:taylor}, and \ref{h:decomposition},
 then Lemma \ref{l:decomposition} implies that
 $\mathcal{L}_{\mu,\ep}^{-1}$ has Fourier symbol
 \begin{align*}
  \hat{\mathcal{L}}^{-1}_{\ep,\mu}(\rho)  = & \frac{1 + \ep^2 \rho^2}{  \beta + \rho^2 [ \ep^2 \beta -\hat{M}(\ep \rho) ] },
  \end{align*} 
  where to simplify notation we have taken $\beta =  \lambda +  (\mu^*+\mu) \rmi n $. 
Because $\hat{ \tilde{K}}_\ep(\xi)$ also satisfies Hypothesis \ref{h:decomposition}, a similar analysis as the one presented in Lemmas \ref{l:invertibleBn}, \ref{l:boundinfinity},  and \ref{l:invertibleL_perp},
shows that for a fixed value of $\ep \neq 0$,
  $\mathcal{L}_{\mu,\ep}^{-1}$ 
is an isomorphisms in $ \mathfrak{H}_{\gamma,n}$.

To show that $\mathcal{L}^{-1}_{0,\mu}:  \mathfrak{H}_{\gamma,n} \rightarrow h^4_{\gamma,n} \oplus b^2_n$
is bounded, we look at the related operator,
\[ \frac{1}{\hat{M}(0)} \mathcal{L}_{0,\mu}   u = \Delta u + \alpha u\qquad \alpha = (\lambda +( \mu^* + \mu)\rmi n) /\hat{M}(0), \]
and the commutative diagram
\begin{center}
\begin{tikzcd}[column sep = large, row sep= huge]
H^{4}_\gamma(\R^2;\C) \arrow[r, "( \Delta + \alpha \rmId )"] \arrow[d, "\langle {\bf x} \rangle^\gamma"] 
&H^2_\gamma(\R^2;\C) \arrow[d, "\langle {\bf x} \rangle^\gamma"] \\
H^4(\R^2;\C) \arrow[r, "A(\gamma)"] 
& H^2(\R^2;\C)
\end{tikzcd}
\end{center}
Notice that the operator
\[ A(\gamma)u = (\Delta +  \alpha \rmId ) u + \left(  2 \nabla u \cdot \nabla \langle x \rangle^{-\gamma} + u \Delta \langle x \rangle^{-\gamma} \right) \langle x \rangle^{\gamma}\]
is a compact perturbation of $(\Delta +  \alpha \rmId ): H^4(\R^2;\C) \longrightarrow H^2(\R^2;\C)$, 
which is invertible since $\alpha \in \C$, with $\mathrm{Im}(\alpha) \neq 0$. It  follows that 
$(\Delta +  \alpha \rmId ): H^4_\gamma(\R^2;\C) \longrightarrow H^2_\gamma(\R^2;\C)$ is Fredholm index zero.
Then, a similar analysis as that in  Lemma B1 of \cite{jaramillo2018} 
shows that this last operator has a trivial kernel, and as a result it is also invertible. 
This proves that $\mathcal{L}^{-1}_{0,\mu}: h^2_{\gamma,n} \rightarrow h^4_{\gamma,n}$ is bounded.

 Similarly, the Fourier symbol for $ \mathcal{L}^{-1}_{0,\mu}(\rho)$, when viewed as a function of 
$\rho $, is analytic and uniformly bounded on strip  in $ \C$ containing the real line. 
Recalling that $\rho = | \xi|$, with $\xi \in \R^2$, it follows that this symbol is also an $L^2(\R^2)$ function,
and that it satisfies the hypotheses of Theorem \ref{t:reed} ( Theorem IX.13 in \cite{reed1975}). As a result $ \hat{ \mathcal{L}}^{-1}_{0,\mu}(\rho)$
generates an $L^1(\R^2)$ convolution kernel, and
by Young's inequality we obtain that
$\mathcal{L}^{-1}_{0,\mu}: W^{2,\infty}(\R^2;\C) \longrightarrow W^{2,\infty}(\R^2;\C)$ is bounded.
Our result
then follows by Lemma \ref{l:diagonal} and the fact that these operators are invariant under rotations.

{\bf Step 2:} Next, we show that the operators are $C^1$ with respect to the parameter $\ep$.

Looking at the symbol  
$\hat{\mathcal{L}}_{\mu,\ep}^{-1}$  one notices that it is
smooth with respect to the parameter $\ep$. 
It then follows that the corresponding operator
is differentiable, and that
 $\partial_\ep \mathcal{L}^{-1}_{\mu, \ep} $ 
has symbol

\begin{align*}
\partial_\ep \hat{\mathcal{L}}^{-1}_{\ep,\mu}(\rho)  = & -\frac{\rho^2}{  \beta + \rho^2 [ \ep^2 \beta - \hat{M}(\ep \rho) ] }
 \left[ -2 \ep 
+ \frac{  ( 1+ \ep^2 \rho^2) }{  \beta + \rho^2 [ \ep^2 \beta - \hat{M}(\ep \rho) ]  } ( 2 \ep \beta - \partial_\ep \hat{M}(\ep \rho) \right].
\end{align*}

To prove that the corresponding operator
\[
\begin{array}{c c c c }
\partial_\ep \mathcal{L}_{\mu,\ep}^{-1} :& \mathfrak{H}_{\gamma,n}& \longrightarrow &   \mathfrak{H}_{\gamma,n},
\end{array}
\]
 is well defined when $\ep \neq 0$,
notice first that  $\hat{M}(\ep \rho)$ and all its derivatives  are uniformly bounded and
 analytic  on a strip containing the real line in the complex plane.
  Therefore, the symbols inside  the brackets satisfy these same properties and thus define isomorphisms on $H^2_\gamma(\R^2)$ for all $\gamma \in \R$.
  Consequently, the main properties of the operator  $\partial_\ep \mathcal{L}_{\mu,\ep}^{-1}$ are determined
by the symbol $\frac{\rho^2}{  \beta + \rho^2 [ \ep^2 \beta - \hat{M}(\ep \rho) ] }$.
It is then straightforward to check that $\partial_\ep \mathcal{L}_{\mu,\ep}^{-1} : H^2_\gamma(\R^2) \longrightarrow H^2_\gamma(\R^2) $ is bounded. Lemma \ref{l:diagonal} then gives the same result for $\partial_\ep \mathcal{L}_{\mu,\ep}^{-1} : h^2_{\gamma,n}(\R^2) \longrightarrow h^2_{\gamma,n}(\R^2) $.

To show that  $\partial_\ep \mathcal{L}_{\mu,\ep}^{-1} : b^2_n(\R^2) \longrightarrow b^2_n(\R^2) $
is bounded, notice that the functions, $ -\frac{\rho^2}{  \beta + \rho^2 [ \ep^2 \beta - \hat{M}(\ep \rho) ] }$ and
$\frac{  ( 1+ \ep^2 \rho^2) }{  \beta + \rho^2 [ \ep^2 \beta - \hat{M}(\ep \rho) ]  } $
have a similar structure as the symbol for the operator $(K\ast + \zeta)^{-1}$ appearing in Lemma \ref{l:boundinfinity}.
As a result they both define isomorphisms on $W^{2,\infty}(\R^2)$. 

A similar result holds for the symbol
$\frac{ \partial_\ep \hat{M}(\ep \rho) ( 1+ \ep^2 \rho^2) }{  \beta + \rho^2 [ \ep^2 \beta - \hat{M}(\ep \rho) ]  } $.
First, Hypothesis \ref{h:decomposition} implies that $\partial_\ep \hat{M}(\ep \rho) \sim \rmO(1/\rho)$
as $\rho \to \infty$. Recalling that $\rho = |\xi|$, as a function of $ \xi \in \R^2$ this symbol defines an $L^2(\R^2)$ function and satisfies the assumptions of Theorem \ref{t:reed}. As in Lemma \ref{l:boundinfinity}, we may then conclude that the corresponding operator is bounded from $W^{2,\infty}(\R^2)$ to $W^{2,\infty}(\R^2)$.
We can therefore view $\partial_\ep \mathcal{L}_{\mu,\ep}^{-1}$ as a composition of 
invertible maps, from $W^{2,\infty}(\R^2)$ back to the same space. Since $b^2_n$ is a closed subset  $W^{2,\infty}(\R^2)$ 
the result then follow.

When $\ep =0$, similar arguments as in the previous two paragraphs show that the operator
$\partial_\ep \mathcal{L}_{\mu,0}^{-1} : \mathfrak{H}_{\gamma,n}  \longrightarrow    h^4_{\gamma,n} \oplus b^2_n$
with symbol
\[ \partial_\ep \mathcal{L}_{\mu,0}^{-1}  = \frac{- \rho^2}{ \beta - \rho^2 \hat{M}(0)}
 \left[   \frac{- \partial_\ep \hat{M}(0)}{  \beta - \rho^2 \hat{M}(0)} \right], \]
is also well defined.

{\bf Step 3:} Finally, we show that the operators are $C^1$ with respect to the parameter $\mu$.

It is clear that $\mathcal{L}_{\mu,\ep}$ is continuously differentiable 
with respect to the parameter $\mu$. 
To show that its inverse has this same property we fix $\ep$ and use the following notation. 
We write $\mathcal{L}_\mu = \mathcal{L}(\mu)$ 
to highlight the dependence of the operator on $\mu$. 
Given $f \in \mathfrak{H}_{\gamma,n}$, 
we let $w(\mu)$ denote the solution to $\mathcal{L}_\mu w =f$ 
and we look at the following equality,
\[ w(\mu + h \mu) - w(\mu) = - \mathcal{L}^{-1}(\mu) \left[  \mathcal{L}(\mu + \mu h)-  \mathcal{L}(\mu)  \right ] w(\mu + h \mu). \]
Since  the operator 
$\left[  \mathcal{L}(\mu + \mu h)-  \mathcal{L}(\mu)  \right ] = \rmi \mu h n$ 
is bounded from $\mathfrak{H}_{\gamma,n}$ back to this same space, 
the above expression then shows that $\mathcal{L}^{-1}_\mu$ is continuous with respect $\mu$. 
At the same time, the above equality shows that the derivative of 
$\mathcal{L}^{-1}(\mu)$ with respect $\mu$ is an operator from 
$\mathfrak{H}_{\gamma,n}$ back to $\mathfrak{H}_{\gamma,n}$,
which is also of the form 
$ - \mathcal{L}^{-1}(\mu) \rmi \mu n \mathcal{L}^{-1}(\mu)$. 
Because this last operator is the composition of maps that depend continuously on $\mu$, 
it is itself also continuous with respect to this parameter.
 \end{proof}

 We finish this section with some comments.
While the constants $\lambda$ and $a$  appearing in the normal form \eqref{e:cgl2}
have specific expressions that can be obtained using a multiple-scale analysis (see Section \ref{s:example} below),
the value of $\mu$ remains an unknown that needs to be determined when solving
this  equation.
As a result, showing existence of specific solutions becomes a nonlinear 
eigenvalue problem, where the form of the solution, $w$, 
has to be determined at the same time as the value of the parameter $\mu$.
A similar problem is encountered when showing the existence of target patterns in 
oscillatory media when an impurity is present.  
There,  the frequency of the waves that emanate from the impurity is not known a priori,
and like $\mu$ here, it is a parameter that needs to be determined and 
thus plays the role of the eigenvalue.
One approach for rigorously finding these target patterns relies on a combination of 
matched asymptotics together with the implicit function theorem, see \cite{jaramillo2018}. 
In this approach, one first shows the existence of target patterns 
for all values of the frequency that lie on a small interval.
Then the matching between the form of the solution in the far field and 
the shape of the solution at intermediate distances provides 
an approximation for the value of frequency that is selected
(which one can show depends on the strength of the impurity). 
We intend to use a similar approach to prove the existence of spiral waves.
The idea would be to first prove the existence of 
solutions to equation \eqref{e:cgl2}, 
valid for all values of $\mu$ on a small interval,
and that have the form $w(r;\mu) = (1 - \rho(r;\mu) ) \rme^{\rmi \phi(r;\mu)}$, with $\partial_r \phi \to kr$ as $r \to \infty$ for some $k \in \R_+$.
Then by matching the form of the solution in the far field with an intermediate approximation 
(that in some sense is unique)
 should allow us to determine the value of $\mu$ that is selected by the system.

 \section{Example}\label{s:example}
In this section we consider the following nonlocal FitzHugh-Nagumo system, posed in $\R^2$,
\begin{align}\label{e:nonlocalFHN}
u_t &= K \ast u + \frac{1}{\tau} ( u - u^3 -v)\\ \nonumber
v_t &= \beta u + \delta.
\end{align}
Here $\tau $ is a small positive parameter, $\beta>0, \delta \in \R$, 
and we assume the convolution operator, $K$, has Fourier symbol 
$$ \hat{K}(\xi) = \dfrac{-\sigma |\xi|^2}{1 + D |\xi|^2}, \quad D,\sigma>0.$$ 

Our goal is to use the methods developed in the previous section, 
together with a multiple scale analysis, to derive a normal form for rotating solutions.

{\bf Set Up:} To start off the multiple scale analysis, 
we first expand the system about the homogeneous steady state 
$(u_*,v_*) = ( \frac{-\delta}{\beta}, (\frac{\delta}{\beta})^3 - \frac{\delta}{\beta} )$, to reveal the linear terms, 
\begin{align*}
u_t &= K \ast u + \frac{1}{\tau} ( (1-  3u_*^2)u -v -3u_*u^2 -u^3),\\
v_t &= \beta u .
\end{align*}
Inserting the rotating solution ansatz, $U(r, \theta,t ) = U(r , \theta + ct)$, 
letting $\lambda = (1 - 3u_*^2)$, 
and writing the resulting equations in matrix form, leads to
\begin{equation}\label{e:FHN}
\begin{bmatrix} c u_\theta \\ c v_\theta \end{bmatrix} = 
\begin{bmatrix}K \ast & 0 \\ 0 & 0 \end{bmatrix}    \begin{bmatrix} u\\ v \end{bmatrix} +
\begin{bmatrix}\lambda \tau^{-1}& -\tau^{-1} \\ \beta & 0 \end{bmatrix}  \begin{bmatrix} u\\ v \end{bmatrix}
+\begin{bmatrix}-\tau^{-1}( 3u_*u^2 + u^3) \\ 0 \end{bmatrix}.
\end{equation}

Using the notation from the previous section, we split the linear terms as
\[
\begin{array}{c c c c c}
A & = & A_0 & + & \ep^2 \bar{\lambda} A_1\\[2ex]
  \begin{bmatrix}\lambda \tau^{-1}& -\tau^{-1} \\ \beta & 0 \end{bmatrix}  &=&
   \begin{bmatrix}0& -\tau^{-1} \\ \beta & 0 \end{bmatrix}  &+& \ep^2 \begin{bmatrix}\bar{\lambda} \tau^{-1}& 0\\ 0 & 0 \end{bmatrix}  .
\end{array}
\]
Then, the eigenvalues of $A_0$ are $\nu = \pm \rmi \sqrt{\beta/\tau}$, 
with corresponding right and left eigenvectors,
\begin{equation}\label{e:vectors}
 W_{1,2} = \begin{bmatrix} - \tau^{-1} \\ \pm \rmi \sqrt{\beta/\tau} \end{bmatrix} \qquad W^*_{1,2} = \frac{1}{2}\begin{bmatrix} -\tau & \mp \rmi \sqrt{\tau/\beta} \end{bmatrix},
 \end{equation}
which are normalized to guarantee that their inner product, 
$\langle W^*_{1,2} , W_{1,2} \rangle =1$. 

 Notice that the nonlinear terms do not depend on the parameter 
 $\lambda$ and can be written as $\tilde{F}(U) = MUU + NUUU$ with
 \begin{align}\label{e:quadratic}
 MUV &= \langle -3\tau^{-1}u_* u_1 v_1, 0 \rangle,\\ \label{e:cubic}
 N UVW & = \langle- \tau^{-1} u_1v_1w_1, 0 \rangle,
 \end{align}
 where $U,V,$ and $W$ are generic vector functions such that 
 $U=(u_1,u_2), V=(v_1,v_2)$ and $W=(w_1,w_2)$.

As in the general case presented in Section \ref{s:normalform}, we let 
$c = c^* + \mu$, where $c^*n_0  = \sqrt{\beta/\tau}$ and $\mu$ is a small parameter. 
Here we are interested in one-armed spirals, so we consider the case when $n_0=1$. 
Again we assume the scalings, 
$s = \ep r$, $\lambda = \ep^2 \bar{\lambda}$, $\mu = \ep^2 \bar{\mu}$  
and consider the expansion
\[  U(r,s) = \ep U_1(s) + \ep^2 U_2(s) + \ep^3 U_3(r), \]
where $U_1 \in X_\parallel$ and $U_{2,3} \in D_\perp$. 
After inserting this ansatz into equation  \eqref{e:FHN} 
and equating coefficients of different powers of $\ep$, one finds that
\begin{align*}
U_1(s) = &W_1 w(s) \rme^{\rmi n_0 \theta} + \overline{W}_1 \overline{w}(s) \rme^{-\rmi n_0 \theta},\\
 U_2(s) =& V_1 w^2 \rme^{2\rmi n_0 \theta} + V_0 |w|^2 + V_{-1} \overline{w}^2 \rme^{-2\rmi n_0 \theta},
 \end{align*}
where the vectors $V_1,V_0,V_{-1}$, satisfy the equations
\begin{align*}
(2\rmi c^*n_0 - A_0) V_1= M W_1W_1 &= -\frac{3 u_*}{\tau^3} \begin{bmatrix} 1 \\0 \end{bmatrix},\\[2ex]
(-2\rmi c^*n_0 - A_0) V_{-1}=  M \overline{W}_1\overline{W}_1& = -\frac{3 u_*}{\tau^3} \begin{bmatrix} 1 \\0 \end{bmatrix},\\[2ex]
- A_0 V_0 =  2M W_1 \overline{W}_1 &= -\frac{6 u_*}{\tau^3} \begin{bmatrix} 1 \\0 \end{bmatrix}.
\end{align*}
They are therefore given by
  \[ V_1 = \frac{u_*}{\tau^2} \begin{bmatrix} 2 \rmi / \sqrt{\beta \tau} \\ 1\end{bmatrix},
  \qquad V_0 =  \frac{u_*}{\tau^2} \begin{bmatrix} -6 \\ 0\end{bmatrix},
  \qquad V_{-1} = \frac{u_*}{\tau^2} \begin{bmatrix} -2 \rmi / \sqrt{\beta \tau} \\ 1\end{bmatrix}.
\]

{\bf Normal Form:} The analysis of the previous section then shows that the normal form for this system is
 \[
 \begin{split}
  \tilde{K}_\ep \ast U_1& -\bar{\mu} \partial_\theta U_1 + \bar{\lambda} A_1(\lambda) U_1 +\\
  & P \left[2 M U_1U_2+ N U_1U_1U_1 + \rmO \Big (\ep (|U_1||U_3| + |U_2 + \ep U3|^4) \Big) \right]=0.
  \end{split} \]
Simplifying this equation using the projection $P$, defined as in \eqref{e:projection},
 then leads to
 \[\tilde{K}_\ep \ast w -  \rmi \bar{\mu} w + \frac{\bar{\lambda}}{\tau} w + (a_1+a_2) |w|^2 w + \rmO(\ep |w|^4w),\]
 and its complex conjugate. In particular, we have that 
\begin{align*}
 a_1 = & \langle W^*_{1} , 2 M(W_1V_0 + \overline{W}_1V_1) \rangle =  \frac{6 u_*^2}{\tau^3} \Big ( 3- \frac{\rmi}{\sqrt{\beta \tau}} \Big), \\
 a_2 = & \langle W^*_{1} , 3 N (W_1W_1 \overline{W}_1) \rangle = - \frac{3}{2 \tau^3},
 \end{align*}
 which are found using the expressions of $U_1$ and $U_2$.

 \section{Discussion}\label{s:discussion}
In this paper, we derived a normal form for systems of equations modeling oscillatory media
with nonlocal coupling. 
Because of their nonlocal nature, one is not able to
use standard techniques from spatial dynamics to obtain this amplitude equation.
The method we use in this paper relies instead on a combination of
Lyapunov-Schmidt reduction and a multiple-scales analysis, which is very
similar to the approach taken in the physics literature. Our main contribution has been to
set up the equations in an appropriate Banach space,  which then allowed us to
decompose the linear part of our system into an invertible 
operator and a bounded operator. 
This decomposition is an essential ingredient for carrying out the Lyapunov-Schmidt reduction,
and for arriving at the normal form. 

In our analysis we also showed the existence of solutions to the normal form equation.
Because this equation is precisely the reduced equation 
obtained from the Lyapunov-Schmidt reduction,
by showing existence of solutions to the normal form 
we also obtain existence of solutions to the full system. 
We emphasize that in contrast to other equations that are 
more commonly referred to amplitude equations,
 say for example the complex Ginzburg-Landau equation,
our normal form equation accounts for all terms that are part of the reduced equation.
This includes higher order terms for which we do not have explicit expressions.
The point here is that even without explicit knowledge of these terms,
we are able to show the existence of solutions to the normal form and
to obtain a first order approximation the solutions of the full system.
Moreover, in contrast to the complex Ginzburg-Landau equation,
which is a parabolic equation
and requires additional analysis to prove the validity 
of its approximations,
the validity of our normal form follows easily from the Lyapunov-Schmidt reduction.
That this is the case follows from the fact that
here we are mainly interested in rotating waves, while the full complex Ginzburg-Landau equation 
gives information about the general-time dependent problem.

To obtain the existence of solutions to the reduced equation, we assumed
that the speed of these solutions corresponds to a free parameter. 
More precisely, the rotational speed of solutions appears
 in the normal form as a the parameter $\mu$.
Here we showed that for all values of $\mu$ in a small interval, 
solutions to the reduced equation exist.
However, as already pointed out in Section \ref{s:validity},
some solutions of interest, like for example spiral waves,
 correspond to specific rotating wave solutions whose 
speed is selected by the system.
This means that in order to find these patterns
one has to view the normal form as a nonlinear eigenvalue problem. 
We remark that this is not a feature of the nonlocal character of the equations, and that a similar
result is seen in the case of other oscillatory systems that are well represented by
reaction diffusion systems. Indeed, in \cite{scheel1998} a center manifold is used 
to derive a similar normal form for reaction diffusion systems undergoing a 
Hopf bifurcation. In this reference, spiral wave solutions are shown to exists using
spatial dynamics and singular perturbation methods. 
In particular, it is shown that in the supercritical case there is a family of spiral wave solutions
which is parametrized by $\mu$, but that in addition there is one particular solution
whose speed is selected by the system. 
On the other hand, in the subcritical case the system
always selects the value of speed. 
We suspect that similar results holds as well in the nonlocal case and we plan to 
address this problem using matched asymptotics and the implicit function theorem in future work.

Finally, we point out that the approach presented here works 
only for finding 2-dimensional rotating wave solutions.
However, the general philosophy of using a multiple-scale analysis combined with
a Lyapunov-Schmidt reduction (based on well-chosen Sobolev spaces) extends to higher dimensions.
For example, one could set up the problem of proving existence of scroll waves  in 
3-d reaction diffusion systems as follows. First, recall that scroll waves
can be modeled as a stack of 2-d spiral waves that
rotate about a filament, see for example \cite{keener1988, keener1992}.
In the simplified case of a straight filament that runs along the z-direction one can consider 
an ansatz of the form
\[ U(r, \theta, z ,t) = V\Big(r, \theta -ct + \psi(\ep z, \ep^2 t) \Big) + \ep^2 u(r, \theta, \ep z, \ep^2 t),\]
where $(r, \theta, z) \in \R^3$ represent cylindrical coordinates and 
$V(r, \theta) = V(r, \theta -ct)$ describes a 2-d spiral wave solution.
Inserting this guess into the original system, 
\[ U_t= \Delta U + F(U), \quad x \in \R^3,\]
results in the expression
\[ L u= V_\theta \frac{\partial \psi}{\partial T} - V_{\theta \theta} \left(\frac{\partial \psi}{\partial Z}\right)^2 - V_{ \theta} \frac{\partial^2 \psi}{\partial Z^2} + \ep^2 \frac{\partial u}{\partial T} - \ep^2 \frac{\partial^2 u}{\partial Z^2}
+ F(U) -F(V) - D_UF(V)u
\]
with $Lu = \Delta_2u + D_UF(V)u$ and $T= \ep^2t, Z=\ep z$. 

The strategy would then be to choose a Banach space $\mathcal{X}$ 
and a projection $P: \mathcal{X} \longrightarrow \mathcal{X}_\perp$, 
$(1-P): \mathcal{X} \longrightarrow \mathcal{X}_\parallel$,
 such that the restriction of the operator $L$ to these subspaces satisfies:
$L_\perp: \mathcal{X}_\perp \rightarrow \mathcal{Y}_\perp$ is invertible, 
while $L_\parallel: \mathcal{X}_\parallel \rightarrow \mathcal{Y}_\parallel$ is bounded. 
Then one can apply Lyapunov-Schmidt reduction to derive a reduced equation, as done in this paper.
Because one-armed spiral wave solutions $V(r, \theta) \longrightarrow \rme^{i kr} \rme^{\rmi \theta} $ as $r \rightarrow \infty$, with $k>0$, the term $D_UF(V)$ approaches a periodic function of $r$. 
Then the decomposition into Fourier modes presented here 
together with the results from \cite{jaramillo2019}, where the Fredholm properties of elliptic
operators with periodic coefficients is established, suggest that such a space, $\mathcal{X}$,
 and projection $P$ might be available.
Notice that by assuming $u \in \mathcal{X}_\perp$, one would recover with this approach
 a Burger-type equation for the modulation $\psi(Z,T)$. That is,
\[ V_\theta \frac{\partial \psi}{\partial T} - V_{\theta \theta} \left(\frac{\partial \psi}{\partial Z}\right)^2 - V_{ \theta} \frac{\partial^2 \psi}{\partial Z^2}  + \rmO(\ep^2) =0.\]

\section{Appendix}

\begin{lemma}
The Fourier Transform maps the spaces  
$$h_n = \{ f \in L^2(\R^2) \mid f( z) = g(r) \rme^{\rmi n \theta} , g \in L^2_r(\R^2)\}$$ 
back to themselves. In particular, given 
$f( z)= f(r \rme^{\rmi \theta}) = g(r) \rme^{\rmi n \theta} \in h_n$, 
then the Fourier transform of these functions can be written as
\[ \mathcal{F}[f(z)] = \mathcal{P}_n[g](\rho) \rme^{\rmi n \phi} = \breve{g}(\rho) \rme^{\rmi n \phi},\]
where
\[ \mathcal{P}_n[g] (\rho) =(-\rmi)^{n} \int_0^\infty g(r) J_n(r \rho) r \;dr, \]
and $J_n(z)$ is the $n$-th order Bessel function of the first kind.
Moreover,
\[ \mathcal{F}^{-1}[\hat{f}(w)] =  \mathcal{P}^{-1}_n[\breve{g}](r) \rme^{\rmi n \theta} = g(r) \rme^{\rmi n \theta},\]
with
\[ \mathcal{P}^{-1}_n[\breve{g}] (r) = \rmi^{n} \int_0^\infty \breve{g}(\rho) J_n(r \rho) \rho \;d\rho. \]
\end{lemma} 

 \begin{proof}
First, we notice that because the Fourier Transform,
$\mathcal{F}$, commutes with orthogonal transformations, 
if $f$ is a radial function then so is $\hat{f}$, 
so that $h_0$ maps back to itself under $\mathcal{F}$. 

Next, given 
$f \in h_n \cap L^1(\R^2)$, i.e. $f(z) = e^{in\theta}g(r)$, 
we want to show that 
$\hat{f}(\rho \rme^{\rmi \phi}) = \rme^{\rmi n \phi} \tilde{f}(\rho)$, 
for some radial $\tilde{f}$. To see why this holds, let $\psi$ be constant and define 
$G(z) = f(r \rme^{\rmi (\theta + \psi)}) $. Then,
$G(z)  = \rme^{\rmi n (\theta + \psi)} g(r) =  \rme^{\rmi n  \psi} f(z) $. 
Therefore, 
\[ \mathcal{F}[G(z)] = \mathcal{F}[ \rme^{\rmi n  \psi} f(z) ]= \rme^{\rmi n  \psi} \hat{f}(w).\] 
On the other hand, because $\rme^{\rmi \psi}$ represents a rotation, 
and the Fourier Transform commutes with orthogonal transformations,
\[ \mathcal{F}[G(z)] = \mathcal{F}[f( \rme^{\rmi \psi} z )] = \hat{f}( \rme^{\rmi \psi} w).\]
This implies that 
$\hat{f}( \rme^{\rmi \psi} w) =  \rme^{\rmi n  \psi} \hat{f}(w)$ 
for all $w$ and all $\psi$.  Letting $w= \rho$ we obtain the desired result for those 
$f \in h_n \cap L^1(\R^2)$. 
Since $h_n \cap L^1(\R^2)$ is dense in $h_n$, 
we can conclude that $\mathcal{F}$ maps the spaces  $h_n$ back to themselves.

Finally, given $f(z) \in h_n $ we have that
\begin{align*}
\mathcal{F}[f(r \rme^{\rmi \theta}) ]  = & \mathcal{F}[g(r) \rme^{\rmi n \theta} ]  \\
= & \frac{1}{2\pi} \int_{\R^2} g(r) \rme^{\rmi n \theta} \rme^{- \rmi \xi \cdot x} \;dx\\
=& \frac{1}{2\pi} \int_0^\infty \int_0^{2\pi} g(r) \rme^{\rmi n \theta} \rme^{- \rmi r \rho \cos( \theta- \phi)} \;d\theta \; r\;dr\\
=& \frac{\rme^{\rmi n \phi}}{2\pi} \int_0^\infty \int_0^{2\pi} g(r) \rme^{\rmi n (\theta-\phi)} \rme^{- \rmi r \rho \cos( \theta- \phi)} \;d\theta \; r\;dr \\
=& \frac{\rme^{\rmi n \phi}}{2\pi} \int_0^\infty g(r) \int_{-\phi+\pi}^{2\pi- \phi +\pi}  \rme^{\rmi n (\psi-\pi)} \rme^{- \rmi r \rho \cos( \psi - \pi)} \;d\psi \; r\;dr 
\end{align*}
where this last integral follows form the change of variables 
$\psi = \theta - \phi + \pi$. If we now focus on the inner integral, 
we notice that because the integrand is $2\pi$-periodic, then
\begin{align*}
\int_{-\phi+\pi}^{2\pi- \phi +\pi}  \rme^{\rmi n (\psi-\pi)} \rme^{- \rmi r \rho \cos( \psi - \pi)} \;d\psi
& = \int_{0}^{2\pi}  \rme^{\rmi n (\psi-\pi)} \rme^{- \rmi r \rho \cos( \psi - \pi)} \;d\psi\\
&= \int_{0}^{2\pi} (-1)^n  \rme^{\rmi n \psi} \rme^{ \rmi r \rho \cos(\psi)} \;d\psi\\
&= \int_{0}^{2\pi} (-1)^n \cos(n \psi) \rme^{ \rmi r \rho \cos(\psi)} \;d\psi\\
& = 2 \pi (-i)^n J_n(\rho r).
\end{align*}
Where in the last line we used the following definition for 
the $n$-th order Bessel function \cite{DLMF}[Eq. 10.9.2]	
\[ J_n(z) =\frac{ (\rmi)^{-n}}{\pi} \int_0^\pi \rme^{\rmi z \cos \theta} \cos(n \theta) \;d\theta.\]
Going back,
\begin{align*}
\mathcal{F}[f(r \rme^{\rmi \theta}) ] =  & \mathcal{F}[g(r) \rme^{\rmi n \theta} ]  \\
= & \frac{\rme^{\rmi n \phi}}{2\pi} \int_0^\infty 2 \pi (-i)^n g(r) J_n(\rho r)\; r\;dr \\
= & \rme^{\rmi n \phi} (-\rmi)^n \int_0^\infty  g(r) J_n(\rho r)\; r\;dr \\
= & \rme^{\rmi n \phi} \mathcal{P}_n[g](\rho)
\end{align*}
A similar calculation then shows that
\[ \mathcal{P}_n^{-1}[\breve{g}] = (\rmi)^n \int_0^\infty \breve{g}(\rho) J_n(r \rho) \rho\; d\rho.\]
That the transformations $\mathcal{P}_n$ and $\mathcal{P}^{-1}_n$ 
are inverses of each other follows from the identity
\[ \int_0^\infty x J_\alpha(ux) J_\alpha(vx) \;dx = \frac{1}{u} \delta(u-v)\]
which holds for $\alpha> -1/2$, see for example \cite{arfken1999}[Sec. 11.2].
\end{proof}

\begin{lemma}
Suppose $f \in M^{2,2}_\gamma(\R^d)$ then $|f(x)| \leq C \|f \|_{M^{2,2}_\gamma} |x|^{-(\gamma+d/2)}$, with $C$ a generic constant.
\end{lemma}
\begin{proof}
Let $(r, \theta) \in (\R^+, \Sigma)$ denote a point in $\R^d$ in spherical coordinates. Given any $f \in M^{2,2}_\gamma(\R^d)$, we may find an upper bound for the $L^2(\Sigma)$ norm of  the function $f(\cdot, R)$, where $R \in \R^+$ is a fixed number, as follows.
\begin{align*}
\| f(\cdot,R)\|^2_{L^2(\Sigma)} = & \int_\Sigma |f(\theta, R)|^2 \;d\theta\\
\leq & \int_\Sigma \left( \int^\infty_R |\partial_rf(\theta,r)| \;dr \right)^2 \; d\theta\\
\leq & \int_\Sigma \left( \int^\infty_R   r^\alpha r^{\gamma +1} |\partial_rf(\theta,r)|\; r^{(d-1)/2} \;dr   \right)^2 \;d\theta\\
\leq & \int_\Sigma \left( \int^\infty_R   r^{2 \alpha} \;dr \right) \left(  \int^\infty_R   r^{2(\gamma +1)} |\partial_rf(\theta,r)|^2 \; r^{(d-1)} \;dr   \right) \;d\theta\\
\leq & C R^{2\alpha +1}
 \int_\Sigma \left(  \int_0^\infty   r^{2(\gamma +1)} |\partial_rf(\theta,r)|^2 \; r^{(d-1)} \;dr   \right) \;d\theta\\
\end{align*}
Where $C$ is a generic constant, we assumed $2\alpha + 1= -(2 \gamma +d) <0$, and we used Cauchy-Schwarz inequality on the fourth line. If instead we had that  $2\alpha + 1= -(2 \gamma +2) >0$, then the above argument can be again carried out, but now the integration in the $r$ variable would be from zero to $R$. This shows that
\[ \| f(\cdot,R)\|_{L^2(\Sigma)}\leq C R^{-(\gamma +d/2)} \| \nabla f \|_{L^2_{\gamma+1}(\R^d)}.\]

On can also repeat the above argument to show that all $\theta$ derivatives satisfy
\[ \| D_\theta f(\cdot,R)\|_{L^2(\Sigma)} \leq C R^{-(\gamma +1 +d/2)} \| D^2 f\|_{L^2_{\gamma+2}(\R^d)}.\]

 We now recall Theorem 5.9 in Adam's book \cite{adams2003} which shows that given $p>1$ and $mp>(d-1)$, and 
$1 \leq q \leq p$, then there is a constant $C = C(m, d,p,q)$ such that
\[\| f(\cdot, R)\|_{L^\infty(\Sigma)} \leq C
 \| f(\cdot, R)\|^\theta_{W^{m,p} (\Sigma)}\; \|  f(\cdot,R)\|^{1-\theta}_{L^q(\Sigma)}\]
with $\theta = (d-1)p/ [ (d-1)p +(mp-(d-1)) q]$. In our case, $m=1$ and $q=p=2$ leading to
\begin{align*}
\| f(\cdot, R)\|_{L^\infty(\Sigma)} \leq & C \|f(\cdot,R)\|^{1-\theta}_{L^2(\Sigma)}\left[ \|f(\cdot,R)\|_{L^2(\Sigma)} + \sum_{|\beta|=1}  \| D^\beta_\theta f(\cdot,R)\|_{L^2(\Sigma)} \right]^\theta\\
\leq & C R^{-(\gamma +d/2)}\| f\|_{M^{2,2}_\gamma(\R^d)}.
\end{align*}

\end{proof}



\bibliographystyle{plain}
\bibliography{spirals}









%
%
%

\medskip
Received xxxx 20xx; revised xxxx 20xx; early access xxxx 20xx.
\medskip

\end{document}